\documentclass[12pt]{article}
\usepackage{amsmath,amscd,amssymb,theorem,enumerate,epsfig,psfrag}
\usepackage{latexsym}
\makeatletter
\def\proclaim #1. #2\par{\medbreak
  \noindent{\bf#1.\enspace}{\sl#2\par}%
  \ifdim\lastskip<\medskipamount \removelastskip\penalty55
        \medskip\fi}
\makeatother
\newcommand{\1}{{{\mathchoice {\rm 1\mskip-4mu l} {\rm 1\mskip-4mu l}
{\rm 1\mskip-4.5mu l} {\rm 1\mskip-5mu l}}}}
\newcommand{\QED}{\hfill$\Box$}

\newtheorem{theorem}{Theorem}[section]
\newtheorem{corollary}[theorem]{Corollary}
\newtheorem{lemma}[theorem]{Lemma}
\newtheorem{proposition}[theorem]{Proposition}

{  
\theorembodyfont{\rmfamily} 
\newtheorem{remark}[theorem]{Remark}
   
}

{
\theorembodyfont{\rmfamily}
}

\def\p{\partial}

\def\b{\beta}
\def\a{\alpha}

\def\eps{\varepsilon}

\def\i{\iota}

\def\o{\omega}
\def\l{\lambda}

\def\L{\Lambda}

\def\S{\Sigma}

\def\cB{{\cal B}}

\def\cE{{\cal E}}
\def\cF{{\cal F}}
\def\cG{{\cal G}}
\def\cH{{\cal H}}

\def\cJ{{\cal J}}
\def\cL{{\cal L}}

\def\cM{{\cal M}}

\def\cP{{\cal P}}
\def\cT{{\cal T}}
\def\cW{{\cal W}}
\def\cX{{\cal X}}

\def\Map{{\rm Map}}

\def\R{{\mathbb R}}

\def\T{{\mathbb T}}
\def\N{{\mathbb N}}
\def\C{{\mathbb C}}

\def\Z{{\mathbb Z}}

\def\area{{\rm area}}

\def\id{{\rm id}}
\def\im{{\rm im\,}}
\def\Om{\Omega}

\newcommand{\inner}[2]{\langle #1, #2\rangle}

\newcommand{{\Emb}}{{\rm Emb}}
\newcommand{{\Diff}}{{\rm Diff}}
\newcommand{{\Vect}}{{\rm Vect}}
\newcommand{{\Cinf}}{{{\mathcal C}^\infty}}
\newcommand{{\Jreg}}{{\mathcal J}_{\rm reg}}
\newcommand{{\Gr}}{{\rm Gr}}
\newcommand{{\ev}}{{\rm ev}}
\newcommand{{\PL}}{{\rm PL}}
\newcommand{{\PU}}{{\rm PU}}
\newcommand{{\U}}{{\rm U}}
\def\Nabla#1{\nabla\kern-.5ex{}_{#1}}
\def\Tilde{\widetilde}

\newcommand{\IMP}{\Longrightarrow}

\newcommand{\INTO}{\hookrightarrow}

\begin{document}

\title{Loops of Lagrangian submanifolds\\
       and pseudoholomorphic discs}

\author{Meike Akveld\;\;\;\;\;\; Dietmar Salamon \\
        ETH-Z\"urich}

\date{21 March 2000}

\maketitle


\begin{abstract}
We examine three invariants of exact loops 
of Lagrangian submanifolds that are modelled 
on invariants introduced by Polterovich
for loops of Hamiltonian symplectomorphisms.  
One of these is the minimal Hofer length 
in a given Hamiltonian isotopy class.  
We determine the exact values of these
invariants for loops of projective Lagrangian planes.
The proof uses the Gromov invariants of 
an associated symplectic fibration
over the $2$-disc with a Lagrangian subbundle
over the boundary. 
\end{abstract}


\section{Introduction}

In this paper we study the Hofer geometry 
for exact loops of Lagrangian submanifolds 
of a symplectic manifold $(M,\o)$.  
Think of such a loop as a submanifold
$
     \L\subset S^1\times M
$
such that the projection $\L\to S^1$ is a submersion
and 
$$
     \L_t := \{z\in M\,|\,(e^{2\pi it},z)\in\L\}
$$
is a Lagrangian submanifold of $M$ for every~$t$. 
The loop is called exact if there exists a 
Hamiltonian isotopy $\psi_t$ of $M$ such that
$\psi_t(\L_0)=\L_t$ for every~$t$. 
The Hofer length of an exact Lagrangian loop 
$\L$ is defined by 
$$
     \ell(\L) 
     := \int_0^1\left(\max_{\L_t}H_t - \min_{\L_t}H_t\right)\,dt,
$$
where the Hamiltonian functions $H_t:M\to\R$ 
are chosen such that the corresponding Hamiltonian 
isotopy $\psi_t:M\to M$ satisfies $\psi_t(\L_0)=\L_t$.
It is interesting to minimize the Hofer length
over the Hamiltonian isotopy class of $\L$.
This infimum will be denoted by 
$$
      \nu(\L) = \nu(\L;M,\o) := \inf_{\L\sim\L'}\ell(\L').
$$

As an explicit example consider the space 
$\cL=\cL(\C P^n,\R P^n)$ of Lagrangian submanifolds
of $\C P^n$ that are diffeomorphic to $\R P^n$. 
It contains the finite dimensional manifold 
$\PL(n+1)$ of projective Lagrangian planes.
The space $\PL(n+1)$ is the orbit of $\R P^n$ 
under the action of $\PU(n+1)$ and its fundamental
group is isomorphic to $\Z_{n+1}$. 
Consider the loop $\L^k\subset S^1\times\C P^n$ 
defined by
\begin{equation}\label{eq:Lk}
     \L^k := \bigcup_{t\in\R}
     \{e^{2\pi it}\}\times\phi_{kt}(\R P^n),
\end{equation}
where 
$
     \phi_t([z_0:\cdots:z_n])
     := [e^{\pi it}z_0:z_1:\cdots:z_n]
$
and $k\in\Z$.
The loops $\L^j$ and $\L^k$ are 
homotopic in $\PL(n+1)$ (as based loops)
if and only if they are Hamiltonian 
isotopic (as free loops) 
if and only if $k-j$ is divisible by $n+1$. 
If $k-j$ is not divisible by $n+1$ then $\L^j$ 
and $\L^k$ can be distinguished by the Maslov index.
More precisely, every Lagrangian loop
$\L\subset S^1\times\C P^n$, with fibres
$\L_t$ Lagrangian isotopic to $\R P^n$, 
has a well defined Maslov index 
$\mu(\L)\in\Z_{n+1}$.
It is defined as the Maslov index of
a smooth map 
$
     u:D=\{z\in\Z\,|\,|z|\le1\}\to M
$
such that $u(e^{2\pi it})\in\L_t$.
Such maps $u$ always exist and the Maslov 
indices of any two such maps differ by an integer
multiple of $n+1$.  It turns out that
\begin{equation}\label{eq:mu-Lk}
     \mu(\L^k) \equiv k\mbox{ mod }n+1.
\end{equation}

In the case $n=1$ the loop $\L^1$ is obtained by rotating a
great circle on the 2-sphere through 180 degrees
around an axis that passes through the circle.
The result is an embedding of the 
Klein bottle into $S^1\times S^2$. 
The image of this embedding is a Lagrangian 
submanifold of $D\times S^2$
with respect to a suitable symplectic form.
In contrast $\Lambda^0$ is a Lagrangian torus in
$D\times S^2$.  
In general, the cases where $n$ is even and where
$n$ is odd are topologically different. 
If $n$ is even, then $\L^k$ is diffeomorphic 
to $S^1\times\R P^n$ for every $k$.
If $n$ is odd then $\L^j$ 
is diffeomorphic to $\L^k$ if and only if
$k-j$ is even, and $\L^k$ is orientable 
if and only if $k$ is even. 
In particular, $\L^k$ is diffeomorphic to 
$\L^0=S^1\times \R P^n$ whenever $k$ is even. 

Fix $k\in\{1,\dots,n\}$ and consider the 
exact Lagrangian loop
$$
     \L := \bigcup_{t\in\R}
     \{e^{2\pi it}\}\times\psi_t(\R P^n),
$$
where 
$$
     \psi_t([z_0:\cdots:z_n])
     := ([z_0:e^{\pi it}z_1:\cdots:e^{\pi it}z_k:z_{k+1}:\cdots:z_n]).
$$
This loop is Hamiltonian isotopic to $\L^k$ 
and it has Hofer length $1/2$,
whereas $\L^k$ has Hofer length $k/2$.
The next theorem asserts that $\L$ minimizes
the Hofer length in its Hamiltonian isotopy class
and hence is a geodesic for the Hofer metric. 

\medskip
\noindent{\bf Theorem~A\,}
{\it Let $\o\in\Om^2(\C P^n)$ denote the Fubini-Study form
that satisfies the normalization condition
$
      \int_{\C P^n}\o^n = 1.
$
Then 
$$
      \nu(\L^k;\C P^n,\o) = \frac12
$$
for $k=1,\dots,n$ and $\nu(\L^0)=0$. 
}

\medskip
\noindent
This is a Lagrangian analogue of a theorem by 
Polterovich~\cite{P1} about loops of Hamiltonian
symplectomorphisms of complex projective space. 
Following~\cite{P1} we introduce two other invariants
of exact Lagrangian loops $\L\subset S^1\times M$ 
that can be expressed in terms of Hamiltonian connection 
$2$-forms $\tau$ on the trivial bundle $D\times M$ that vanish 
over $\L$.  Let $\cT(\L)\subset\Om^2(D\times M)$ 
denote the space of such connection $2$-forms.
The {\bf relative K-area} $\chi(\L)$ is obtained 
by minimizing the Hofer norm of the curvature $\Om_\tau$
over $\cT(\L)$. The third invariant is related to the 
relative cohomology classes $[\tau]\in H^2(D\times M,\L;\Z)$
of $\tau\in\cT(\L)$.  These form a $1$-dimensional affine 
space parallel to the subspace generated by the 
integral cohomology class $\sigma := [dx\wedge dy/\pi]$.
For $\tau_0,\tau_1\in\cT(\L)$ define $s(\tau_1,\tau_0)\in\R$
by 
$
     s(\tau_1,\tau_0)\sigma = [\tau_1]-[\tau_0].
$
The invariant $\eps(\L)$ is defined by 
$$
     \eps(\L) := \eps^+(\tau_0,\L) - \eps^-(\tau_0,\L),
$$
for $\tau_0\in\cT(\L)$, where 
$$
     \eps^+(\tau_0,\L) 
     := \inf\{s(\tau,\tau_0)\,|\,
       \tau\in\cT(\L),\,\tau^{n+1}>0\},
$$
$$
     \eps^-(\tau_0,\L) 
     := \sup\{s(\tau,\tau_0)\,|\,
       \tau\in\cT(\L),\,\tau^{n+1}<0\}.
$$

\medskip
\noindent{\bf Theorem~B\,}
{\it For every exact Lagrangian loop $\L\subset S^1\times M$}
$$
     \eps(\L)\le\chi(\L)=\nu(\L).
$$
\noindent
A lower bound for $\eps(\L)$ can sometimes be 
obtained by studying pseudoholomorphic sections of 
$D\times M$ with boundary values in $\L$. 
We assume that the pair $(M,\L_0)$ is monotone
and fix a class $A\in H_2(D\times M,\L;\Z)$
that satisfies
$$
     n\pm\mu_\L(A) \le N-2,
$$
where $n=\dim\,\L_0=\dim\,M/2$, $N$ denotes
the minimal Maslov number of the pair $(M,\L_0)$,
and $\mu_\L$ denotes the Maslov class. 
Under these assumptions we define Gromov invariants 
$$
     \Gr^\pm_A(\L) \in H_{n\pm\mu_\L(A)}(\L_0;\Z_2).
$$
A connection $2$-form $\tau\in\cT(\L)$ and 
an $\o$-compatible almost complex structure $J$
on $M$ determine an almost complex structure
$\tilde J=\tilde J(\tau,J)$ on $D\times M$.  
Under our assumptions the moduli space of 
$\tilde J(\tau,\pm J)$-holomorphic
sections of $D\times M$ is, for a generic $\tau$,
a compact smooth manifold of dimension $n\pm\mu_\L(A)$.
The Gromov invariant is defined as the image of the 
mod-2 fundamental class under the evaluation map
$u\mapsto u(1)$. Now let $\L^k\subset S^1\times\C P^n$ 
be given by~(\ref{eq:Lk}) with $1\le k\le n$. 
Let $A^\pm\in H_2(D\times\C P^n,\L^k;\Z)$ be 
the homology classes of the constant sections 
$u^+(x,y)\equiv[1:0:\cdots:0]$ and 
$u^-(x,y)\equiv[0:\cdots:0:1]$.

\medskip
\noindent{\bf Theorem~C\,}
{\it $\Gr^\pm_{A^\pm}(\L^k)\ne 0$. 
}

\medskip
\noindent
Theorem~C can be interpreted as an existence result 
for pseudoholomorphic sections and we shall use this
to prove that $\eps(\L^k)\ge1/2$.
On the other hand the Hamiltonian isotopy class 
of $\L^k$ contains a loop of length equal to $1/2$.  
Hence Theorem~A follows from Theorem~B. 

We expect that the same techniques can be used to obtain 
similar results for general symplectic quotients of 
$\C^n$ by subgroups of ${\rm U}(n)$. These quotients will not, 
in general, satisfy our assumption of monotonicity
for the definition of the Gromov invariants.
However, it should be possible to derive the same conclusions 
by using the invariants introduced in
Cieliebak--Gaio--Salamon~\cite{CGS} instead.
This programme will be carried out elsewhere.

In~\cite{P1,P2,P3,P4} Polterovich studied the Hofer length
of loops 
$
     \psi_t=\psi_{t+1}:M\to M
$ 
of Hamiltonian symplectomorphisms.
Let $P\to S^2$ denote the Hamiltonian fibration
associated to the Hamiltonian loop. 
Poltervich introduced invariants
$\nu^\pm(P)$, $\chi^\pm(P)$, and $\eps^\pm(P)$
on which our invariants are modelled.
Here $\nu^+(P)$ is obtained by 
minimizing the positive part of the 
Hofer length in a given Hamiltonian
isotopy class, the K-area $\chi^+(P)$ 
is a symplectic analogue of an invariant
introduced by Gromov~\cite{G2}, 
and the invariant $\eps^+(P)$
is based on the coupling construction of
Guillemin--Lerman--Sternberg~\cite{GLS}. 
In~\cite{P1,P2} Polterovich proves that 
these invariants are equal:
$$
      \eps^\pm(P) = \chi^\pm(P) = \nu^\pm(P).
$$
We adopt the convention $\pm\nu^\pm(P)\ge 0$. 
Let us denote by $\nu(P)$, $\chi(P)$, and $\eps(P)$ 
the Hamiltonian analogues of our invariants 
of Lagrangian loops. These were also considered 
by Polterovich and he noted that 
$$
      \eps(P) 
      = \eps^+(P)-\eps^-(P)
      = \nu^+(P)-\nu^-(P)
      \le \nu(P).
$$
This is the Hamiltonian analogue of Theorem~B. 
Now consider the Lagrangian loop 
$\L\subset S^1\times\bar M\times M$ given by
$$
      \L_t = {\rm graph}(\psi_t).
$$
The invariants introduced by Polterovich
are related to our invariants by
$$
      \nu(\L)\le\nu(P),\qquad
      \eps(\L)\le\eps(P).
$$
The Gromov invariants of the fibration $P$ associated
to a Hamiltonian loop were independently studied by
Seidel~\cite{S2,S3,S4} and his results were used 
by Lalonde--McDuff--Polterovich~\cite{LMP} to prove 
that Hamiltonian loops act trivially on homology. 
Our results on the Gromov invariants
can be viewed as Lagrangian analogues of results
in~\cite{P1,S2} on the Gromov invariants 
of symplectic fibrations. 

The present paper is organized as follows.  
In Section~\ref{sec:lag} we discuss background material 
about the Hofer metric. The space of Lagrangian 
submanifolds is naturally foliated by Hamiltonian 
isotopy classes and the Hofer metric 
is defined on each leaf of this foliation. 
In Section~\ref{sec:inv} we introduce the invariants
$\nu(\L)$, $\chi(\L)$, and $\eps(\L)$ of exact
Lagrangian loops and give a proof of Theorem~B. 
In the $2$-dimensional case the invariant $\nu(\L)$ 
can sometimes be computed explicitly.  This is done 
in Section~\ref{sec:torus} for the $2$-torus.
In Section~\ref{sec:gro} we introduce the Gromov invariants
and in Section~\ref{sec:cx} we prove Theorems~A and~C.
In Appendix~\ref{app:2} we prove a
result about Hamiltonian isotopy on Riemann surfaces
which is used in Section~\ref{sec:torus}. 

\medskip
\noindent {\bf Acknowledgement:} 
We would like to thank Leonid Polterovich
for suggesting the topic and for many helpful 
discussions.  


\section{The Hofer metric for Lagrangian submanifolds}\label{sec:lag}

Let $(M,\o)$ be a $2n$-dimensional symplectic
manifold and $L$ be a compact connected $n$-manifold
without boundary.  Denote by 
$$
     \cX = \left\{\i\in\Emb(L,M)\,|\,\i^*\o=0\right\}
$$
the space of Lagrangian embeddings of $L$ into $M$. 
The group $\cG=\Diff(L)$ acts on this space by 
$\i\mapsto\i\circ\phi$ for $\phi\in\cG$.  
Two Lagrangian embeddings $\i_0,\i_1\in\cX$
lie in the same $\cG$-orbit if and only if 
they have the same image $\Lambda=\i_0(L)=\i_1(L)$.
Hence the quotient space
$$
     \cL := \cX/\cG
$$
can be naturally identified with the set of Lagrangian 
submanifolds of $M$ that are diffeomorphic to $L$.
A function $\R\to\cL:t\mapsto\Lambda_t$ 
is called {\bf smooth} if there exists a smooth 
function $\R\times L\to M:(t,q)\mapsto\i_t(q)$
such that $\i_t(L)=\Lambda_t$ for all $t$.
One can think of $\cL$ as an infinite dimensional manifold.  

\begin{lemma}\label{le:tangent}
The tangent space of $\cL$ at a point $\Lambda\in\cL$
can be naturally identified with the space 
of closed $1$-forms on $\Lambda$:
$$
     T_\Lambda\cL 
     = \left\{\beta\in\Om^1(\Lambda)\,|\,
       d\beta=0\right\}
$$
\end{lemma}

\begin{proof}
Let $\R\times L\to M:(t,q)\mapsto\i_t(q)$ 
be a smooth function such that $\i_t\in\cX$
for all $t$ and define 
\begin{equation}\label{eq:alpha}
     \alpha_t := \o(v_t,d\i_t\cdot)\in\Om^1(L),\qquad
     v_t:=\p_t\i_t\in\Cinf(L,{\i_t}^*TM).
\end{equation}
Then 
$$
     0 = \p_t{\i_t}^*\o = d\alpha_t
$$
and hence the tangent space of $\cX$ at $\i$ is given by
$$
     T_\i\cX 
     = \left\{v\in\Cinf(L,\i^*TM)\,|\,\o(v,d\i\cdot)\in\Om^1(L) 
        \mbox{ is closed}\right\}.
$$
The tangent space to the $\cG$-orbit consists of all
vector fields of the form $v=d\i\circ\xi$, where 
$\xi\in\Vect(L)$.  The map $v\mapsto\o(v,d\i\cdot)$ 
identifies the quotient space $T_\i\cX/T_\i(\i\cdot\cG)$
with the space of closed $1$-forms on~$L$. 

If $\i_t,\i_t'\in\cX$ are two smooth
paths in $\cX$ that satisfy
$
    \i_t' = \i_t\circ\phi_t
$
for some path $\phi_t\in\cG$ then the vector fields
$
    v_t:=\p_t\i_t
$
and
$ 
    v_t':=\p_t\i_t'
$
are related by
$$
    v_t' = v_t\circ\phi_t + d\i_t\circ\xi_t\circ\phi_t
$$
where $\xi_t\in\Vect(L)$ generates the diffeomorphism
$\phi_t$ via 
$
    \p_t\phi_t = \xi_t\circ\phi_t.
$
Hence the $1$-forms 
$
     \alpha_t := \o(v_t,d\i_t\cdot)
$
and
$
     \alpha_t' := \o(v_t',d\i_t'\cdot)
$
are related by 
$$
    \alpha_t' = {\phi_t}^*\alpha_t.
$$
Hence two closed $1$-forms $\alpha,\alpha'\in\Om^1(L)$ 
corresponding to two Lagrangian embeddings 
$\i$ and $\i'=\i\circ\phi$ represent the same 
tangent vector of $\cL$ if and only if $\alpha'=\phi^*\alpha$
or, equivalently, $\i_*\alpha={\i'}_*\alpha'$.
This proves the lemma. 
\QED\medskip
\end{proof}

Let $\R\to\cL:t\mapsto\Lambda_t$ be a smooth path of
Lagrangian submanifolds.  We define the derivative of this
path at time $t$ by
$$
     \p_t\Lambda_t := {\i_t}_*\alpha_t,
$$
where the path $\R\to\cX:t\mapsto\i_t$ is chosen such that
$\i_t(L)=\Lambda_t$ for every $t$ and $\alpha_t$
is defined by~(\ref{eq:alpha}). The proof of 
Lemma~\ref{le:tangent} shows that the 
$1$-form $\beta_t={\i_t}_*\alpha_t\in\Om^1(\Lambda_t)$
is closed and is independent of the choice of the lift
$t\mapsto\i_t$ used to define it. 

We wish to study Hamiltonian isotopies of Lagrangian submanifolds.
This corresponds to paths in $\cL$ that are tangent to the 
subbundle
$$
    \cH = \left\{(\Lambda,\beta)\in T\cL\,|\,
    \Lambda\in\cL,\,\beta\in\Om^1(\Lambda)\mbox{  is exact}
    \right\}.
$$
Abstractly, one can think of $\cH$ as a distribution
on $\cL$.  It follows from Weinstein's Lagrangian 
neighbourhood theorem that this distribution is integrable. 
We shall see that the leaf through $\Lambda_0\in\cL$ 
consists of all Lagrangian submanifolds 
of $M$ that are Hamiltonian isotopic to $\Lambda_0$ .  
To be more precise, let $\R\times M\to\R:(t,z)\mapsto H_t(z)$ 
be a smooth Hamiltonian function and denote by 
$\R\times M\to M:(t,z)\mapsto\psi_t(z)$ the 
Hamiltonian isotopy generated by $H$ via
\begin{equation}\label{eq:ham}
    \frac{d}{dt}\psi_t = X_t\circ\psi_t,\qquad
    \i(X_t)\o = dH_t,\qquad
    \psi_0=\id.
\end{equation}

\begin{lemma}\label{le:ham1}
Let $\R\to\cL:t\mapsto\L_t$ be a smooth path of Lagrangian
submanifolds and $\psi_t$ be a Hamiltonian isotopy on $M$
generated by the Hamiltonian functions $H_t:M\to\R$
via~(\ref{eq:ham}). Then 
$\Lambda_t = \psi_t(\Lambda_0)$ for every $t$
if and only if 
$$
     \p_t\Lambda_t = dH_t|_{\Lambda_t}
$$
for every $t$. 
\end{lemma}

\begin{proof}
Choose a smooth path $\R\to\cX:t\mapsto\i_t$ 
such that $\i_t(L)=\Lambda_t$ for every $t$ and let  
$\alpha_t\in\Om^1(L)$ be defined by~(\ref{eq:alpha}).
Then $\p_t\Lambda_t = dH_t|_{\Lambda_t}$ 
if and only if $d(H_t\circ\i_t)=\alpha_t$.  
It follows from the definitions 
that this is equivalent to 
$$
     X_t(\i_t(q))-\p_t\i_t(q) \in \im d\i_t(q)
$$
for all $t$ and all $q$. This means that there exists
a smooth family of vector fields $\xi_t\in\Vect(L)$
such that 
$$
     X_t\circ\i_t = \p_t\i_t + d\i_t\circ\xi_t.
$$
Equivalently,
$
     \psi_t\circ\i_0 = \i_t\circ\phi_t,
$
where the isotopy $\phi_t\in\Diff(L)$ is generated by 
$\xi_t$ via 
$
     \p_t\phi_t=\xi_t\circ\phi_t
$
and
$
     \phi_0=\id.
$
This proves the lemma.
\QED\medskip
\end{proof}

The previous lemma shows that every path in $\cL$ that is 
generated by a Hamiltonian isotopy is tangent to $\cH$. 
The converse is  proved next. 

\begin{lemma}\label{le:ham2}
A smooth path $[0,1]\to\cL:t\mapsto\Lambda_t$ 
is tangent to $\cH$ if and only if 
there exists a Hamiltonian isotopy $t\mapsto\psi_t$
such that $\psi_t(\Lambda_0)=\Lambda_t$ for every~$t$.
\end{lemma}

\begin{proof}
The ``if'' part was proved in Lemma~\ref{le:ham1}. 
Suppose that the path $t\mapsto\Lambda_t$ is tangent 
to $\cH$. Choose a smooth function 
$[0,1]\to\cX:t\mapsto\i_t$ such that $\i_t(L)=\Lambda_t$
for every $t$ and let  
$\alpha_t\in\Om^1(L)$ be defined by~(\ref{eq:alpha}).
By assumption, $\alpha_t$ is exact for every $t$. 
Fix a smooth path $q_t\in L$ and, for every~$t$, 
choose $h_t:L\to\R$ such that
$$
     dh_t = \alpha_t,\qquad h_t(q_t)=0.
$$
Then the function $\R\times L\to\R:(t,q)\mapsto h_t(q)$
is smooth. We construct a smooth function 
$[0,1]\times M\to\R:(t,z)\mapsto H_t(z)$
such that
\begin{equation}\label{eq:ham1}
     H_t\circ\i_t = h_t.
\end{equation}
Choose an almost complex structure $J$ on $M$
that is compatible with $\o$. 
Let $\eps>0$ be so small that, 
for every $t\in[0,1]$, the map 
$$
     T\Lambda_t\to M:
     (z,v)\mapsto\exp_z(Jv)
$$
restricts to a diffeomorphism from the $\eps$-neighbourhood
of the zero section in $T\Lambda_t$ onto the open neighbourhood
$$
     U_t := \left\{\exp_z(Jv)\,|\,
     z\in\Lambda_t,\,v\in T_z\Lambda_t,\,|v|<\eps\right\}
$$
of $\Lambda_t$ in $M$. Choose a cutoff function 
$\rho:[0,\eps]\to[0,1]$ such that $\rho(r)=1$ 
for $r<\eps/3$ and $\rho(r)=0$ for $r>2\eps/3$.
Define $H_t:M\to\R$ by
$$
     H_t(\exp_z(Jv)) 
     := \rho(|v|)h_t\circ{\i_t}^{-1}(z)
$$
for $z\in\Lambda_t$ and $v\in T_z\Lambda_t$ with $|v|<\eps$, 
and by $H_t(z):=0$ for $z\in M\setminus U_t$.
Then $H_t$ satisfies~(\ref{eq:ham1}) and hence
$$
     dH_t|_{\L_t} = {\i_t}_*dh_t = {\i_t}_*\a_t = \p_t\L_t.
$$
By Lemma~\ref{le:ham1}, the Hamiltonian isotopy $\psi_t$
generated by $H_t$ satisfies 
$
     \psi_t(\Lambda_0)=\Lambda_t
$
for every $t$. This proves the lemma.
\QED
\end{proof}

\begin{remark}\label{rmk:norm}
The Hamiltonian functions constructed in 
Lemma~\ref{le:ham2} satisfy
\begin{equation}\label{eq:norm}
     \max H_t = \max h_t,\qquad
     \min H_t = \min h_t
\end{equation}
for every $t$. With a slightly more sophisticated 
argument one can show that the Hamiltonian functions
can be chosen such that the Hamiltonian vector fields
$X_t$ satisfy $\p_t\i_t = X_t\circ\i_t$ and hence
the resulting Hamiltonian isotopy satisfies
\begin{equation}\label{eq:i-psi}
     \psi_t\circ\i_0=\i_t.
\end{equation}
However, in general there does not exist a Hamiltonian
isotopy that satisfies both~(\ref{eq:norm}) 
and~(\ref{eq:i-psi}).
\end{remark}

\begin{lemma}\label{le:natural}
Let $\R\to\cL:t\mapsto\Lambda_t$ be a smooth
path of Lagrangian submanifolds. 
Let $\R\to\Diff(M,\o):t\mapsto\psi_t$ be a symplectic
isotopy and define $\beta_t\in\Om^1(M)$ by 
$
     \beta_t := \i(Y_t)\o,
$
where
$
     \p_t\psi_t = Y_t\circ\psi_t.
$
Then $\beta_t$ is closed and the path 
$\Lambda_t':={\psi_t}^{-1}(\Lambda_t)$ satisfies
$$
     \p_t\Lambda_t'
     = {\psi_t}^*\left(\p_t\Lambda_t 
     - \beta_t|_{\Lambda_t}\right).
$$
\end{lemma}

\begin{proof}
Choose a lift $\R\to\cX:t\mapsto\i_t$ of $t\mapsto\Lambda_t$
and denote
$$
      \i_t' := {\psi_t}^{-1}\circ\i_t,\qquad
      \a_t  := \o(\p_t\i_t,d\i_t\cdot),\qquad
      \a_t' := \o(\p_t\i_t',d\i_t'\cdot).
$$
Then 
$
      \a_t' = \a_t - {\i_t}^*\beta_t
$
and hence 
$$
      \p_t\Lambda_t'
= 
      {\i_t'}_*\a_t'  
= 
      {\psi_t}^*{\i_t}_*\a_t'  
= 
      {\psi_t}^*\left({\i_t}_*\a_t - \beta_t\right)  
= 
      {\psi_t}^*\left(\p_t\Lambda_t - \beta_t\right)
$$
as claimed. 
\QED\medskip
\end{proof}

\smallbreak

The subbundle $\cH\subset T\cL$ carries a natural norm.  
Following Hofer~\cite{H} we define the norm
of an exact $1$-form $\alpha=dh\in\Om^1(\L)$ by 
$$
    \|dh\| := \max h - \min h.
$$
This norm gives rise to a distance function on each
leaf of the foliation determined by $\cH$.  
Let $\cL_0$ be such a leaf.  By Lemma~\ref{le:ham2},
$\cL_0$ is the Hamiltonian isotopy class of 
any Lagrangian submanifold $\Lambda\in\cL_0$. 
Let $[0,1]\to\cL_0:t\mapsto\Lambda_t$ be a smooth
path in $\cL_0$. The length of this path is defined by 
$$
    \ell(\{\Lambda_t\})
    := \int_0^1 \|\p_t\L_t\|\,dt.   
$$
Lemma~\ref{le:ham2} and Remark~\ref{rmk:norm} show that 
\begin{equation}\label{eq:length1}
    \ell(\{\Lambda_t\})
    = \inf_{\psi_t(\Lambda_0)=\Lambda_t}\ell(\{\psi_t\}),
\end{equation}
where the infimum runs over all Hamiltonian isotopies
$t\mapsto\psi_t$ that satisfy
$\psi_t(\Lambda_0)=\Lambda_t$ for all $t$ and 
$\ell(\{\psi_t\})$ denotes the Hofer length (cf.~\cite{H}). 

Now let $\Lambda,\Lambda'\in\cL_0$ and denote by 
$\cP(\Lambda,\Lambda')$ the space of all smooth
paths $[0,1]\to\cL_0:t\mapsto\Lambda_t$ 
that connect $\Lambda_0=\Lambda$ to $\Lambda_1=\Lambda'$. 
The distance between $\Lambda$ and $\Lambda'$ is defined by 
\begin{equation}\label{eq:d}
    d(\Lambda,\Lambda')
    := \inf_{\{\Lambda_t\}\in\cP(\Lambda,\Lambda')} 
    \ell(\{\Lambda_t\}).
\end{equation}
It follows immediately from~(\ref{eq:length1}) that
\begin{equation}\label{eq:d1}
    d(\Lambda,\Lambda')
    = \inf_{\psi(\Lambda)=\Lambda'}
      d(\id,\psi)
\end{equation}
where the infimum runs over all Hamiltonian symplectomorphisms
$\psi$ of $M$ that satisfy $\psi(\Lambda)=\Lambda'$
and $d(\id,\psi)$ denotes the Hofer distance (cf.~\cite{H}). 
The function~(\ref{eq:d}) is obviously nonnegative, 
symmetric, and satisfies the triangle inequality.
That it defines a metric is a deep theorem due to
Chekanov~\cite{C}. 

\begin{theorem}[Chekanov]
If $\Lambda\ne\Lambda'$ then $d(\Lambda,\Lambda')>0$. 
\end{theorem}



\begin{remark}
In~\cite{MILINKOVICH} Milinkovi\'c studied 
geodesics in the space of Lagrangian submanifolds. 
Generalizing a result by Bialy and Polterovich~\cite{BP},
he proved that the distance of two exact Lagrangian
submanifolds $\L={\rm graph}(dS)$ and $\L'={\rm graph}(dS')$
of the cotangent bundle $T^*L$ is given by
$$
       d(\L,\L') = \|d(S-S')\|.
$$
\end{remark}


\section{Invariants of Lagrangian loops}\label{sec:inv}

In this section we shall consider exact loops of 
Lagrangian submanifolds.  In the terminology of the 
previous section this corresponds to loops inside 
a leaf of the foliation of $\cL$ determined by $\cH$. 
We shall construct three invariants of Hamiltonian 
isotopy classes of such loops and study the relations
between them.


\subsection{The minimal length}

Continue the notation of Section~\ref{sec:lag}.
A {\bf Lagrangian loop} in $M$ is a smooth function
$\R\to\cL:t\mapsto\Lambda_t$ such that 
$$
      \Lambda_{t+1}=\Lambda_t
$$
for all $t\in\R$. Such a loop
determines a subset $\Lambda\subset S^1\times M$
defined by
\begin{equation}\label{eq:lambda}
      \Lambda
      := \left\{(e^{2\pi it},z)\,|\,t\in\R,z\in\Lambda_t\right\}.
\end{equation}
Note that a loop $\R\to\cL:t\mapsto\Lambda_t$ is smooth
if and only if this set $\Lambda$ is a smooth 
submanifold of $S^1\times M$. We shall frequently identify
the loop $\R\to\cL:t\mapsto\Lambda_t$ with the corresponding
submanifold $\Lambda\subset S^1\times M$.

A Lagrangian loop $t\mapsto\Lambda_t$ 
is called {\bf exact} if it is tangent to $\cH$, i.e. 
$\p_t\Lambda_t\in\Om^1(\Lambda_t)$ is exact for every $t$. 
Two exact Lagrangian loops $t\mapsto\Lambda_t$ and 
$t\mapsto\Lambda_t'$ are called {\bf Hamiltonian isotopic}
if there exists a smooth function 
$[0,1]\times\R\to\cL:(s,t)\mapsto\Lambda_{s,t}$
such that
$$
      \Lambda_{0,t} = \Lambda_t,\qquad
      \Lambda_{1,t} = \Lambda_t',
$$
the map $t\mapsto\Lambda_{s,t}$ 
is an exact Lagrangian loop for every $s$, 
and $\p_s\Lambda_{s,t}\in\Om^1(\Lambda_{s,t})$ 
is exact for all $s$ and $t$.  Here the function 
$[0,1]\times\R\to\cL:(s,t)\mapsto\Lambda_{s,t}$
is called smooth if there exists a smooth function
$[0,1]\times\R\times L\to M:(s,t,q)\mapsto\i_{s,t}(q)$
such that $\i_{s,t}(L)=\L_{s,t}$ for all 
$s$ and $t$.  Let $\L,\L'\subset S^1\times M$ 
be two exact Lagrangian loops.
We write $\L\sim\L'$ iff $\L$ is Hamiltonian 
isotopic to $\L'$. A Hamiltonian isotopy class corresponds 
to a component in the free loop space of a leaf $\cL_0\subset\cL$
of the foliation determined by $\cH$. 
To every such Hamiltonian isotopy class 
we assign the real number 
$$
     \nu(\L) := \inf_{\L'\sim\L}\ell(\L').
$$
So $\nu(\L)$ is obtained by minimizing the Hofer length over all
exact Lagrangian loops that are Hamiltonian isotopic to $\L$.


\subsection{The relative K-area}

Following Polterovich~\cite{P2}
we introduce the notion of relative K-area. 
This invariant is defined in terms of Hamiltonian
connections on the symplectic fibre bundle $D\times M\to D$
that preserve the subbundle $\Lambda\subset D\times M$ 
defined by~(\ref{eq:lambda}). 
Here $D\subset\C$ denotes the closed unit disc.
We begin by recalling the basic notions
of symplectic connections and curvature
(cf.~\cite{GLS,MS1}).
Think of a connection on $D\times M$ as a horizontal
distribution. Any such connection is determined by a 
{\bf connection $2$-form} on $D\times M$ of the form
$$
      \tau 
      = \o + \a\wedge dx + \b\wedge dy 
        + f dx \wedge dy
$$
where $\a=\a_{x,y}\in\Om^1(M)$, $\b=\b_{x,y}\in\Omega^1(M)$,
and $f=f_{x,y}\in\Om^0(M)$ depend smoothly on $x+iy\in D$.
The horizontal subspace is the $\tau$-orthogonal complement
of the vertical subspace.  Explicitly, the horizontal lifts
of $\p/\p x$ and $\p/\p y$ at $(x+iy,z)\in D\times M$ are the 
vectors $(1,X_{x,y}(z))$ and $(i,Y_{x,y}(z))$, respectively,
where the vector fields $X=X_{x,y},Y=Y_{x,y}\in\Vect(M)$ 
are defined by 
$$
     \i(X)\o=\alpha,\qquad \i(Y)\o = \beta.
$$
Thus the connection associated to $\tau$ is independent of $f$.  
It is called {\bf symplectic} if $\alpha_{x,y}$ and $\beta_{x,y}$ 
are closed for all $x+iy\in D$, and {\bf Hamiltonian} 
if $\alpha_{x,y}$ and $\beta_{x,y}$ are exact for all $x+iy\in D$
and $\tau$ is closed.\footnote
{
In~\cite{MS1} a connection is called 
Hamiltonian if parallel transport along every {\bf loop}
in the base is a Hamiltonian symplectomorphism.  
In the case of a simply connected base this is equivalent
to the existence of a closed $2$-form $\tau$ that represents
this connection.  In contrast, we call a connection Hamiltonian
if parallel transport along every {\bf path} is a Hamiltonian
symplectomorphism.  This notion only makes sense when the bundle 
in question is equipped with a trivialization. 
}
Thus a Hamiltonian connection $2$-form has the form
\begin{equation}\label{eq:tau}
     \tau = \o + dF\wedge dx + dG\wedge dy 
     + (\p_xG - \p_yF+c) dx\wedge dy,
\end{equation}
where $F,G:D\times M\to\R$ and $c:D\to\R$ are smooth maps
such that the functions $F_{x,y}=F(x+iy,\cdot)$ and 
$G_{x,y}=G(x+iy,\cdot)$ have mean value zero:
$$
     \int_MF_{x,y}\o^n = \int_MG_{x,y}\o^n = 0.
$$
In~(\ref{eq:tau}) the $d$ in $dF$ denotes the differential on $M$, 
i.e. $dF$ denotes the smooth family $x+iy\mapsto dF_{x,y}$
of $1$-forms on $M$, and similarly for $dG$. 
We shall only consider Hamiltonian connections with the 
property that parallel transport along the boundary 
preserves~$\Lambda$. 

\begin{lemma}\label{le:con1}
Let $\tau$ be a Hamiltonian connection $2$-form on 
$D\times M$ of the form~(\ref{eq:tau}) and denote
\begin{equation}\label{eq:Ht}
     H_t
     := - 2\pi\sin(2\pi t)F_{\cos(2\pi t),\sin(2\pi t)}
        + 2\pi\cos(2\pi t)G_{\cos(2\pi t),\sin(2\pi t)}.
\end{equation}
Let $\R\to\cL:t\mapsto\Lambda_t$ be an exact Lagrangian loop,
let $\Lambda\subset D\times M$ be defined by~(\ref{eq:lambda}),
and choose a smooth function $\i:\R\times L\to M$ such that 
$\i_t(L)=\Lambda_t$, where $\i_t:=\i(t,\cdot)$.
Then the following are equivalent. 
\begin{description}
\item[(i)]
Parallel transport of $\tau$ along the boundary preserves $\Lambda$. 
\item[(ii)]
$\i^*\tau=0$. 
\item[(iii)]
$
     dH_t|_{\L_t} = \p_t\L_t
$
for every $t\in\R$.
\end{description}
\end{lemma}

\begin{proof}
The parallel transport of $\tau$ along a curve 
$
     t\mapsto x(t)+iy(t)
$
is determined by the 
Hamiltonian functions
$$
     H_t = \dot x(t) F_{x(t),y(t)} + \dot y(t)G_{x(t),y(t)}
$$
via~(\ref{eq:ham}).  The functions $H_t$ in~(\ref{eq:Ht}) 
correspond to the path $t\mapsto e^{2\pi it}$.  
By Lemma~\ref{le:ham1}, the Hamiltonian isotopy determined
by $H_t$ preserves $\L$ if and only if 
$dH_t|_{\L_t}=\p_t\L_t$ for every $t$.
This shows that~(i) is equivalent to~(iii). 

To prove the equivalence of~(ii) and~(iii) note that
$$
     \tau(\p_t\i_t,d\i_t\cdot)
     = \o(\p_t\i_t-X_t\circ\i_t,d\i_t\cdot),
$$
where $X_t\in\Vect(M)$ denotes the Hamiltonian vector
field of $H_t$ as in~(\ref{eq:ham}). 
The right hand side vanishes if and only if 
$
     dH_t|_{\L_t} = \p_t\L_t
$
and the left hand side vanishes if and only if 
$\i^*\tau=0$.  This proves the lemma.
\QED\medskip
\end{proof}

For every exact Lagrangian loop $\R\to\cL:t\mapsto\Lambda_t$
let us denote the set of Hamiltonian connections 
that preserve $\Lambda$ by  
$$
    \cT(\Lambda)
    = \left\{\tau\in\Om^2(D\times M)\,|\,
      \tau\mbox{ has the form }(\ref{eq:tau}),\,
      \tau|_{T\Lambda}=0
      \right\}.
$$
We shall prove in Lemma~\ref{le:con2} below that this set
is nonempty. Let $\R\to\cL:t\mapsto\Lambda_t'$ be another
exact Lagrangian loop. 
A diffeomorphism 
$$
     \Psi:(D\times M,\L)\to(D\times M,\L')
$$ 
is called a {\bf fibrewise (Hamiltonian) symplectomorphism} 
if it has the form
$
     \Psi(x+iy,z) = (x+iy,\psi_{x,y}(z)),
$
where $\psi_{x,y}:M\to M$ is a (Hamiltonian) 
symplectomorphism for all $x,y$. 
In the case $\Lambda=\Lambda'$ we denote by 
$\cG(\Lambda)$ the group of fibrewise Hamiltonian
symplectomorphisms of $(D\times M,\L)$.  
This group acts on $\cT(\Lambda)$ by $\tau\mapsto\Psi^*\tau$.
The {\bf curvature} of a connection $2$-form
$\tau$ of the form~(\ref{eq:tau}) is the function 
$\Om_\tau:D\times M\to\R$ defined by 
\begin{equation}\label{eq:curvature}
    \Om_\tau(x,y,z)
    := \{F_{x,y},G_{x,y}\}(z) + \p_yF_{x,y}(z) - \p_xG_{x,y}(z)
\end{equation}
for $x+iy\in D$ and $z\in M$. It is sometimes useful 
to think of the curvature as a $2$-form 
$\Om_\tau dx\wedge dy$ on $D\times M$ rather than a function.

\begin{lemma}\label{le:con2}
\begin{description}
\item[(i)]
For every exact Lagrangian loop 
$\R\to\cL:t\mapsto\L_t$ the set
$\cT(\L)$ is nonempty.
\item[(ii)]
Two exact Lagrangian loops $\L$ and $\L'$ 
are Hamiltonian isotopic if and only if the corresponding 
pairs $(D\times M,\Lambda)$ and $(D\times M,\Lambda')$
are fibrewise Hamiltonian symplectomorphic.
\item[(iii)]
If $\tau$ is a Hamiltonian connection $2$-form
on $D\times M$ and $\Psi$ is a fibrewise Hamiltonian 
symplectomorphism of $D\times M$ then
$$
     \Om_{\Psi^*\tau} = \Om_\tau\circ\Psi.
$$
\end{description}
\end{lemma}

\begin{proof} Let $\phi_t\in\Diff(L)$ be defined by 
$$
     \i_{t+1}\circ\phi_t = \i_t.
$$
Since $L$ is connected there exists a smooth path
$\R\to L:t\mapsto q_t$ such that, for every $t\in\R$,
\begin{equation}\label{eq:q}
     q_{t+1} = \phi_t(q_t).
\end{equation}
For example choose $q_t$ in the interval 
$0\le t\le 1$ such that $q_t=q_1$ for $1-\eps\le t\le1$ and 
$q_t={\phi_t}^{-1}(q_1)$ for $0\le t\le\eps$.  Then define 
$q_t$ for $t\in\R$ such that~(\ref{eq:q}) is satisfied.
Let $h_t:L\to\R$ be defined by 
$$
     dh_t = \alpha_t := \o(\p_t\i_t,d\i_t\cdot),\qquad
     h_t(q_t) = 0.
$$
By~(\ref{eq:q}), the function $t\mapsto\i_t(q_t)$ 
is $1$-periodic in $t$ and the proof of Lemma~\ref{le:tangent}
shows that the $1$-forms ${\i_t}_*\alpha_t$ are $1$-periodic in~$t$.
Hence the functions $h_t\circ{\i_t}^{-1}$ are $1$-periodic in $t$
and hence, so are the functions $H_t$ defined in
the proof of Lemma~\ref{le:ham2}.  
Now define
$$
     \Tilde{H}_t(z) 
     := H_t(z) - \frac{\int_M H_t\o^n}{\int_M\o^n}.
$$
Let $\rho:[0,1]\to[0,1]$ be a smooth cutoff
function such that $\rho(r)=0$ for $r<\eps$ and $\rho(r)=1$
for $r>1-\eps$ and define $\tau$ by 
\begin{equation}\label{eq:tau-r}
      \Phi^*\tau = \o + \rho(r)d\Tilde{H}_t dt 
      + \dot\rho(r)\Tilde{H}_t dr\wedge dt,
\end{equation}
where $\Phi:[0,1]\times[0,1]\times M\to D\times M$ 
is given by 
$
     \Phi(r,t,z) = (re^{2\pi it},z).
$
Explicitly, $\tau$ has the form~(\ref{eq:tau}) 
where $F,G:D\times M\to\R$ are given by
\begin{equation}\label{eq:FG}
     F_{x,y} = \frac{-\sin(2\pi t)\rho(r)}{2\pi r}\Tilde{H}_t,\qquad
     G_{x,y} = \frac{\cos(2\pi t)\rho(r)}{2\pi r}\Tilde{H}_t,
\end{equation}
for $x+iy=re^{2\pi it}$.
These functions have mean value zero and 
satisfy~(\ref{eq:Ht}) with $H_t$ replaced by~$\Tilde{H}_t$. 
Since $H_t\circ\i_t=h_t$ it follows as in the proof of 
Lemma~\ref{le:ham2} that
$$
     d\Tilde{H}_t|_{\L_t}
     = dH_t|_{\L_t}
     =\p_t\L_t.
$$
By Lemma~\ref{le:con1}, the parallel 
transport of $\tau$ along the boundary preserves $\Lambda$.
Hence $\tau$ is an element of $\cT(\Lambda)$.  
This proves~(i). 

We prove~(ii).  Assume first that there exists a
fibrewise Hamiltonian symplectomorphism of the form
$
      \Psi(x+iy,z)=(x+iy,\psi_{x+iy}(z))
$
such that 
$$
      \psi_{e^{2\pi it}}(\Lambda_t) = \Lambda_t'
$$
for every $t$.   Define 
$$
      \psi_{s,t} := \psi_{se^{2\pi it}},\qquad
      \Lambda_{s,t} := \psi_{s,t}(\Lambda_t)
$$
for $0\le s\le 1$ and $t\in\R$. Then $t\mapsto\Lambda_{s,t}$
is an exact Lagrangian loop for every $s$ and
$\p_s\Lambda_{s,t}\in\Om^1(\Lambda_{s,t})$ is exact 
for all $s$ and $t$.  Hence the Lagrangian loop
$\Lambda_{1,t}=\Lambda_t'$ is Hamiltonian isotopic to 
$\Lambda_{0,t}=\psi_0(\Lambda_t)$. Since $\psi_0$ 
is a Hamiltonian symplectomorphism, the loop $t\mapsto\psi_0(\Lambda_t)$
is Hamiltonian isotopic to $t\mapsto\Lambda_t$.
Conversely, suppose that $t\mapsto\Lambda_t$ and 
$t\mapsto\Lambda_t'$ are two exact Lagrangian loops that are
Hamiltonian isotopic.  Choose an exact isotopy
$(s,t)\mapsto\Lambda_{s,t}$ such that 
$
      \Lambda_{0,t}=\Lambda_t,
$
$
      \Lambda_{1,t}=\Lambda_t',
$
and $\p_s\Lambda_{s,t}=0$ for $s\le1/2$. 
As in the proof of~(i), 
one can construct a smooth family of 
Hamiltonian functions $H_{s,t}:M\to\R$ such that 
$$
      H_{s,t+1}=H_{s,t},\qquad
      dH_{s,t}|_{\Lambda_{s,t}} = \p_s\Lambda_{s,t}.
$$
Define the Hamiltonian symplectomorphisms $\psi_{s,t}:M\to M$ 
by 
$$
      \p_s\psi_{s,t} = X_{s,t}\circ\psi_{s,t},\qquad
      \i(X_{s,t})\o = dH_{s,t},\qquad
      \psi_{0,t}=\id.
$$
Then $\psi_{s,t}=\id$ for $s\le1/2$ and the required
fibrewise Hamiltonian symplectomorphism is given by
$\Psi(se^{2\pi it},z):=(se^{2\pi it},\psi_{s,t}(z))$. 

We prove~(iii). Let $\tau$ be given by~(\ref{eq:tau})
and suppose that 
$$
    \Psi(x+iy,z)=(x+iy,\psi_{x,y}(z))
$$
is a fibrewise Hamiltonian symplectomorphism.
Choose smooth functions $A,B:D\times M\to\R$ such that
the functions
$
     A_{x,y}:=A(x+iy,\cdot)
$
and 
$
     B_{x,y}:=B(x+iy,\cdot)
$
have mean value zero and the Hamiltonian vector fields
$
     X_A=X_{A_{x,y}}
$
and
$
     X_B=X_{B_{x,y}}
$
satisfy
\begin{equation}\label{eq:AB}
    \p_x\psi = X_A\circ\psi,\qquad
    \p_y\psi = X_B\circ\psi.
\end{equation}
Then 
$$
    \Psi^*\tau 
    = \o + d\Tilde{F}\wedge dx + d\Tilde{G}\wedge dy
      + (\p_x\Tilde{G}-\p_y\Tilde{F}+c)dx\wedge dy,
$$
where 
$$
    \Tilde{F} = (F-A)\circ\Psi,\qquad
    \Tilde{G} = (G-B)\circ\Psi.
$$
Hence 
\begin{eqnarray*}
    \Om_{\Psi^*\tau}
&= &
    \p_x\Tilde{G} - \p_y\Tilde{F} - \{\Tilde{F},\Tilde{G}\}  \\
&= &
    \p_x(G-B)\circ\Psi + d(G-B)\circ X_A\circ\Psi  \\
&&
    - \p_y(F-A)\circ\Psi - d(F-A)\circ X_B\circ\Psi \\
&&
    - \{(F-A),(G-B)\}\circ\Psi \\
&= &
    (\p_xG - \p_yF - \{F,G\})\circ\Psi \\
&&
    - (\p_xB-\p_yA-\{A,B\})\circ\Psi  \\
&= &
    \Om_\tau\circ\Psi.
\end{eqnarray*}
The last equality follows from the definition of $A$
and $B$ in~(\ref{eq:AB}). This proves the lemma. 
\QED\medskip
\end{proof}

The {\bf relative K-area} of an exact 
Lagrangian loop $\Lambda$ is defined by 
$$
    \chi(\Lambda) 
    := \inf_{\tau\in\cT(\Lambda)}\left\|\Om_\tau\right\|,
$$
where
$$ 
    \|\Om_\tau\| 
    := \int_D \left(\max_{z\in M}\Om_\tau(x,y,z)
       - \min_{z\in M}\Om_\tau(x,y,z)\right)\,dxdy.
$$

\begin{theorem}\label{thm:nu-chi}
For every exact Lagrangian loop 
$\L\subset S^1\times M$
$$
     \chi(\L) = \nu(\L).
$$
\end{theorem}

\begin{proof}
Let $\R\to\cL:t\mapsto\Lambda_t$ be an exact Lagrangian 
loop.  Let $\tau\in\Om^2(D\times M)$ be the connection $2$-form 
defined by~(\ref{eq:tau-r}) in the proof of Lemma~\ref{le:con2}, 
where the cutoff function $\rho:[0,1]\to[0,1]$ is chosen 
to be nondecreasing.  Then 
$$
     \Phi^*(Fdx+Gdy) = \rho\Tilde{H}_t dt,
$$
where $F,G:D\times M\to\R$ are given by~(\ref{eq:FG})
and $\Phi(r,t,z)=(re^{2\pi it},z)$.
Taking the differential of this $1$-form on $[0,1]^2\times M$
we find
$$
     \Phi^*((\p_xG-\p_yF)dx\wedge dy)
     = \dot\rho\Tilde{H}_t dr\wedge dt.
$$
Since $\{F,G\}=0$ and $\Phi^*(dx\wedge dy)=2\pi r dr\wedge dt$
we obtain
$$
     \Om_\tau(re^{2\pi it},z)
     = -\frac{\dot\rho(r)}{2\pi r}\Tilde{H}_t(z).
$$
Moreover,
$$
     \left\|\Tilde{H}_t\right\|
     = \max_M\Tilde{H}_t - \min_M\Tilde{H}_t
     = \max_{\Lambda_t}\Tilde{H}_t - \min_{\Lambda_t}\Tilde{H}_t,
$$
and hence
$$
     \left\|\Om_\tau\right\| 
     = \int_0^1\int_0^1 
       \dot\rho(r)\left\|\Tilde{H}_t\right\|
       \,drdt  
     = \int_0^1\left\|\Tilde{H}_t\right\|\,dt  
     = \ell(\L).
$$
This implies
$
     \chi(\L) \le \ell(\L).
$
If $\L$ and $\L'$ are Hamiltonian isotopic
then, by Lemma~\ref{le:con2}~(ii), there exists a 
fibrewise Hamiltonian symplectomorphism 
$\Psi$ of $D\times M$ such that $\Psi(\L)=\L'$.
Hence $\tau\in\cT(\L')$ if and only if 
$\Psi^*\tau\in\cT(\L)$.  By Lemma~\ref{le:con2}~(iii),
$
     \chi(\L)=\chi(\L')\le\ell(\Lambda').
$
Hence
$
     \chi(\L) \le \nu(\L).
$

We prove that $\nu(\L)\le\chi(\L)$. 
Let $\tau\in\cT(\L)$.
We shall construct an exact Lagrangian 
loop $\L'$ that is Hamiltonian
isotopic to $\L$ and satisfies
\begin{equation}\label{eq:L'}
     \ell(\L') \le \|\Om_\tau\|.
\end{equation}
Suppose that $\tau$ has the form~(\ref{eq:tau}).
Since the function $c$ in~(\ref{eq:tau}) 
has no effect on the curvature we may assume,
without loss of generality, that $c\equiv0$. 
Define $H=H_{r,t}:M\to\R$ and $K=K_{r,t}:M\to\R$ 
by the formula 
$$
    \Phi^*\tau = \o + dK\wedge dr + dH\wedge dt
    + (\p_rH-\p_tK)dr\wedge dt.
$$
Explicitly,
$$
    K_{r,t} = \cos(2\pi t)F_{re^{2\pi it}} 
    + \sin(2\pi t)G_{re^{2\pi it}},
$$
$$
    H_{r,t} = 2\pi r\cos(2\pi t)G_{re^{2\pi it}} 
    - 2\pi r\sin(2\pi t)F_{re^{2\pi it}}.
$$
Define the Hamiltonian symplectomorphisms
$\psi_{r,t}:M\to M$ by 
$$
    \p_r\psi_{r,t} = X_{K_{r,t}}\circ\psi_{r,t},\qquad
    \psi_{0,t} = \id.
$$
Then the loop 
$$
    \Lambda_t' = {\psi_{1,t}}^{-1}(\Lambda_t)
$$
is evidently Hamiltonian isotopic to $\Lambda$.
We shall prove that it satisfies~(\ref{eq:L'}). 
To see this, denote by $\Psi$ the fibrewise Hamiltonian
symplectomorphism of $[0,1]^2\times M$ given by
$$
    \Psi(r,t,z) = (r,t,\psi_{r,t}(z)).
$$
Then, as in the proof of Lemma~\ref{le:con2},
we obtain 
$$
    \Psi^*\Phi^*\tau 
    = \o + dH'\wedge dt + \p_rH'dr\wedge dt,
$$
where $H'_{r,t}=(H_{r,t}-B_{r,t})\circ\psi_{r,t}$
and $B_{r,t}:M\to\R$ is defined by 
$\p_t\psi_{r,t}=X_{B_{r,t}}\circ\psi_{r,t}$.
These functions satisfy 
$$
    \|\Om_\tau\| = \int_0^1\int_0^1
    \|\p_rH'_{r,t}\|\,drdt,\qquad
    H'_{0,t} = 0.
$$
Moreover, by Lemma~\ref{le:natural}, we have
\begin{eqnarray*}
    \p_t\Lambda_t'
&= &
    {\psi_{1,t}}^*\left(\p_t\Lambda_t 
    - dB_{1,t}|_{\Lambda_t}\right) \\
&= &
    {\psi_{1,t}}^*(dH_{1,t}|_{\Lambda_t})
    - d(B_{1,t}\circ\psi_{1,t})|_{\Lambda_t'}  \\
&= &
    dH_{1,t}'|_{\Lambda_t'}.
\end{eqnarray*}
Hence the length of $\Lambda'$ is given by
\begin{eqnarray*}
     \ell(\L') 
&= &
     \int_0^1 
     \left(\max_{\Lambda_t'}H_{1,t}'
     - \min_{\Lambda_t'}H_{1,t}'\right)
     \,dt   \\
&\le &
     \int_0^1 
     \left(\max_MH_{1,t}'
     - \min_MH_{1,t}'\right)
     \,dt   \\
&= &
     \int_0^1 
     \left(
     \max_M\left(\int_0^1\p_rH'_{r,t}\,dr\right)
     - \min_M\left(\int_0^1\p_rH'_{r,t}\,dr\right)
     \right)
     \,dt   \\
&\le &
     \int_0^1\int_0^1 
     \left(
     \max_M\p_rH'_{r,t}
     - \min_M\p_rH'_{r,t}
     \right)
     \,drdt   \\
&= &
     \int_0^1\int_0^1 
     \left\|\p_rH'_{r,t}\right\|
     \,drdt   \\
&= &
     \|\Om_\tau\|.
\end{eqnarray*}
Thus we have proved that for every $\tau\in\cT(\L)$ 
there exists an exact Lagrangian loop $\L'$ that 
is Hamiltonian isotopic to $\L$ and satisfies
$\ell(\L')\le\|\Om_\tau\|$.  Hence $\chi(\L)\le\nu(\L)$
and this proves the theorem.
\QED
\end{proof}


\subsection{The non-symplectic interval}\label{sec:nonsymp}
 
Let $\Lambda\subset D\times M$ be an exact Lagrangian loop
and $\tau\in\cT(\Lambda)$ be a Hamiltonian connection
$2$-form. Since $\tau$ is closed and vanishes on $\Lambda$
(see Lemma~\ref{le:con1}) it determines a relative 
cohomology class
$$
      [\tau]\in H^2(D\times M,\Lambda;\R).
$$
Let $\Sigma$ be a compact oriented Riemann surface
with (possibly empty) boundary $\p\Sigma$. 
A smooth map $v:(\Sigma,\p\Sigma)\to (D\times M,\L)$ 
determines a $2$-dimensional relative homology
class
$$
      [v] := v_*[\Sigma]\in H_2(D\times M,\L;\Z).
$$
The pairing of this class with $[\tau]$ is given by
$$
      \inner{[\tau]}{[v]} = \int_\Sigma v^*\tau.
$$
Since every 2-dimensional integral homology class 
of the pair $(D\times M,\L)$ can be represented by a 
smooth map $v$ as above, the cohomology class $[\tau]$
is uniquely determined by these pairings. 
Define $\sigma\in H^2(D\times M,\Lambda;\R)$ by
\begin{equation}\label{eq:sigma}
     \inner{\sigma}{[v]} = \deg(\pi\circ v)
\end{equation}
for every $v:(\Sigma,\p\Sigma)\to(D\times M,\L)$,
where 
$$
     \pi:(D\times M,\Lambda)\to(D,\p D)
$$
denotes the obvious projection. 
In~(\ref{eq:sigma}) the degree of a smooth map
$v_0:(\S,\p\S)\to(D,\p D)$ is understood as the 
degree of its restriction to the boundary. 
It agrees with the number of preimages of an 
interior regular value, counted with appropriate 
signs (cf. Milnor~\cite{MILNOR}). Note that 
$$
     \sigma = \frac{1}{\pi}[dx\wedge dy]
$$
and hence $\sigma$ agrees with the pullback of the 
positive integral generator of $H^2(D,\p D;\R)$ 
under the projection $\pi$. 

\begin{lemma}\label{le:coh}
Let $\tau_0,\tau_1\in\cT(\L)$.  Then there exists 
a constant $s=s(\tau_1,\tau_0)\in\R$ such that 
$$
     [\tau_1]-[\tau_0] = s\sigma.
$$
\end{lemma}

\begin{proof}
Let $\tau_i$ be given by~(\ref{eq:tau}) with $F,G,c$ 
replaced by $F_i,G_i,c_i$ for $i=0,1$. Denote 
$$
    F:=F_1-F_0,\qquad G:=G_1-G_0,\qquad c:=c_1-c_0,
$$
and let $H_t:M\to\R$ be defined by~(\ref{eq:Ht}). 
Since $\tau_0,\tau_1\in\cT(\L)$ it follows from 
Lemma~\ref{le:con1} that there exists a function 
$h:\R/\Z\to\R$ such that 
$$
      H_t|_{\L_t}\equiv h(t)
$$
for every $t\in\R$.  We shall prove that the 
required identity holds with
$$
      s := \int_0^1h(t)\,dt + \int_D c\,dxdy.
$$
To see this note that, by~(\ref{eq:Ht}),
\begin{equation}\label{eq:hdt}
      (Fdx+Gdy)|_\L = \pi^*\a_h.
\end{equation}
where $\alpha_h\in\Om^1(S^1)$ denotes the pushforward
of the $1$-form $hdt\in\Om^1(\R/\Z)$ under the 
diffeomorphisms $\R/\Z\to S^1:[t]\mapsto e^{2\pi it}$.
Let $\S$ be a compact
oriented Riemann surface and $v:\Sigma\to D\times M$ 
be a smooth function such that $v(\p\S)\subset\L$.
Denote $v_0:=\pi\circ v:(\S,\p\S)\to(D,\p D)$. Then
\begin{eqnarray*}
    \int_\S v^*(\tau_1-\tau_0)
&= &
    \int_\S v^*
    \left(dF\wedge dx + dG\wedge dy + (\p_xG-\p_yF+c)dx\wedge dy\right) \\
&= &
    \int_{\p\S} v^*(Fdx+Gdy)
    + \int_\S {v_0}^*(cdx\wedge dy)  \\
&= &
    \int_{\p\S} {v_0}^*\a_h
    + \int_\S {v_0}^*(cdx\wedge dy)  \\
&= &
    s\deg(v_0).
\end{eqnarray*}
The penultimate equality follows from~(\ref{eq:hdt})
and the last from the identities
\begin{equation}\label{eq:deg1}
    \int_{\p\S} {v_0}^*\a_h
    = \deg(v_0)\int_{S^1}\a_h
\end{equation}
and 
\begin{equation}\label{eq:deg2}
    \int_\S {v_0}^*(cdx\wedge dy)
    = \deg(v_0)\int_Dcdx\wedge dy.
\end{equation}
Here~(\ref{eq:deg1}) is the degree theorem
for maps between compact $1$-manifolds 
and~(\ref{eq:deg2}) is the degree theorem 
for maps between $2$-manifolds with boundary.  
More precisely, if the function $c:D\to\R$
has mean value zero then there exists a $1$-form 
$\alpha\in\Om^1(D)$ such that 
$
      d\alpha = cdx\wedge dy
$
and
$
      \alpha|_{TS^1}=0.
$
This implies that the left hand side of~(\ref{eq:deg2})
vanishes. Hence it suffices to establish~(\ref{eq:deg2})
for constant functions $c$ and this reduces
to~(\ref{eq:deg1}).  This proves the lemma. 
\QED\medskip
\end{proof}

Let $\tau_0\in\cT(\L)$. We shall now address the question
which cohomology classes $[\tau_0]+s\sigma$ can be represented
by nondegenerate Hamiltonian connection $2$-forms.
Such a $2$-form is a symplectic form on $D\times M$ 
with respect to which $\Lambda$ is a Lagrangian submanifold.
Denote
$$
     \cT^\pm(\L)
     := \left\{\tau\in\cT(\L)\,|\,\pm\tau^{n+1}>0\right\}.
$$
Here the inequality $\tau^{n+1}>0$ means that 
$\tau^{n+1}=f\,dx\wedge dy\wedge \o^n$,
where $f:D\times M\to\R$ is a positive function.
For $\tau_0\in\cT(\L)$ we define
$$
     \eps^+(\tau_0,\L)
     := \inf\left\{s(\tau,\tau_0)\,|\,\tau\in\cT^+(\L)\right\},
$$
$$
     \eps^-(\tau_0,\L)
     := \sup\left\{s(\tau,\tau_0)\,|\,\tau\in\cT^-(\L)\right\}.
$$
The proof of Theorem~\ref{thm:eps-nu} below shows that
the class $[\tau_0]+s\sigma$ can be represented
by a symplectic form $\tau\in\cT^\pm(\L)$
for $\pm s$ sufficiently large
and hence $\pm\eps^\pm(\tau_0,\L)<\infty$. 
Evidently, $\eps^\pm(\tau_1,\L)-\eps^\pm(\tau_0,\L)=s(\tau_1,\tau_0)$.
Hence the number
$$
     \eps(\Lambda) := \eps^+(\tau_0,\Lambda) - \eps^-(\tau_0,\Lambda)
$$
is independent of the connection $2$-form $\tau_0\in\cT(\L)$
used to define it. This number is called the
{\bf width of the nonsymplectic interval}.

\begin{theorem}\label{thm:eps-nu}
For every exact Lagrangian loop $\L\subset D\times M$
$$
    \eps(\L) \le \chi(\L).
$$
\end{theorem}

\begin{proof}
Let $\R\to\cL:t\mapsto\L_t$ be an exact Lagrangian
loop and $F,G:D\times M\to\R$ be smooth functions
such that the functions $H_t:M\to\R$ 
defined by~(\ref{eq:Ht}) satisfy
$
     dH_t|_{\L_t} = \p_t\L_t
$
for every $t$. For every smooth function $c:D\to\R$ let
$\tau_c\in\cT(\L)$ be given by~(\ref{eq:tau}). 
In particular, $\tau_0$ is given by~(\ref{eq:tau}) with $c=0$. 
We shall prove that 
\begin{equation}\label{eq:eps+}
     \eps^+(\tau_0,\L)
     \le \int_D\max_{z\in M}\Om_{\tau_0}(x,y,z)\,dxdy,
\end{equation}
\begin{equation}\label{eq:eps-}
     \eps^-(\tau_0,\L)
     \ge \int_D\min_{z\in M}\Om_{\tau_0}(x,y,z)\,dxdy.
\end{equation}
To see this, note that
$$
     ndF\wedge dG\wedge\o^{n-1}
     = \{F,G\}\o^n
$$
and hence
\begin{equation}\label{eq:tau-n+1}
\begin{array}{rcl}
     \tau^{n+1}
&= &
     (n+1)(\p_xG-\p_yF+c)dx\wedge dy\wedge\o^n  \\
&&
     +\, n(n+1)dF\wedge dx\wedge dG\wedge dy\wedge\o^{n-1} \\
&= &
     (n+1)(\p_xG-\p_yF-\{F,G\}+c)dx\wedge dy\wedge\o^n \\
&= &
     (n+1)(c-\Om_{\tau_0})dx\wedge dy\wedge\o^n.
\end{array}
\end{equation}
This shows that $\tau_c$ is nondegenerate if and only if 
$
     c(x,y) \ne \Om_{\tau_0}(x,y,z)
$
for all $(x+iy,z)\in D\times M$. 
Fix a number 
$$
     s > \int_D\max_{z\in M}\Om_{\tau_0}(x,y,z)\,dxdy.
$$
Choose a smooth function $c:D\to\R$ such that 
$$
     c(x,y) > \max_{z\in M}\Om_{\tau_0}(x,y,z)
$$
for all $x+iy\in D$ and 
$$
     \int_Dc\,dxdy = s.
$$
Then $\tau_c$ is nondegenerate
and represents the class 
$
     [\tau_c] = [\tau_0] + s\sigma.  
$
This proves~(\ref{eq:eps+})
and~(\ref{eq:eps-}) follows from a similar argument. 
It follows from~(\ref{eq:eps+}) and~(\ref{eq:eps-})
that
\begin{eqnarray*}
     \eps(\L)
&= &
     \eps^+(\tau_0,\Lambda) - \eps^-(\tau_0,\Lambda) \\
&\le &
     \int_D
     \left(
     \max_{z\in M}\Om_{\tau_0}(x,y,z)
     - \min_{z\in M}\Om_{\tau_0}(x,y,z)
     \right)
     \,dxdy  \\
&= &
     \|\Om_{\tau_0}\|.
\end{eqnarray*}
Since the curvature of $\tau_0$ is equal to the curvature 
of $\tau_c$ for every $c$ it follows that 
$\eps(\L)\le\|\Om_\tau\|$ for every $\tau\in\cT(\L)$
and hence $\eps(\L)\le\chi(\L)$. 
This proves the theorem.
\QED
\end{proof}

\begin{remark}\label{rmk:nonsymplectic}
Let us denote
\begin{equation}\label{eq:T}
    T(\L) :=
    \left\{[\tau]\in H^2(D\times M,\L;\R)\,|\,
    \tau\in\cT(\L)\right\}.
\end{equation}
By Lemma~\ref{le:coh}, this set is a $1$-dimensional
affine subspace of $H^2(D\times M,\L;\R)$.
Denote
$$
    T^\pm(\L) :=
    \left\{[\tau]\,|\,
    \tau\in\cT^\pm(\L)\right\}.
$$
These sets are open and connected.  
To prove connectedness, 
let $\tau_i\in\cT^+(\L)$ be given by~(\ref{eq:tau})
with $F,G,c$ replaced by $F_i,G_i,c_i$ for $i=0,1$.
By~(\ref{eq:tau-n+1}), $c_i>\Om_{\tau_i}$.
Assume without loss of generality that
$s(\tau_1,\tau_0)\ge0$. Then the path
$
     [0,1]\to T^+(\L):
     t\mapsto[\tau_0]+ts(\tau_1,\tau_0)\sigma
$
connects $[\tau_0]$ with $[\tau_1]$. 
This shows that the sets $T^\pm(\L)$ are connected.
The complement $T(\L)\setminus(T^-(\L)\cup T^+(\L))$
is compact and connected. 
It can be expressed in the form
$$
     T(\L)\setminus(T^-(\L)\cup T^+(\L))
     = \left\{[\tau_0]+s\sigma\,|\,
       \eps^-(\tau_0,\L)\le s\le\eps^+(\tau_0,\L)
       \right\}
$$
for every $\tau_0\in\cT(\L)$. 
We do not know if this complement is always
nonempty or, equivalently, if $\eps(\L)$ 
is always nonnegative.  
\end{remark}


\section{Loops on the 2-torus}\label{sec:torus}

\begin{figure}[htp]
\centerline{\psfig{figure=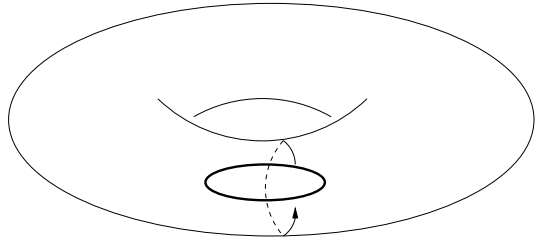}} 
\caption{{A Lagrangian loop on the 2-torus}}\label{fig:torus1}
\end{figure}

Consider the torus $M=\T^2=\R^2/\Z^2$ 
with the standard symplectic form 
$$
     \o = dx\wedge dy
$$
and let $\pi:\R^2\to\T^2$ denote the projection.
Let 
$$
     B_r=\{(s,t)\in\R^2\,|\,s^2+t^2\le r^2\}
$$
and suppose that $S\subset\T^2$ is 
the image of an embedding $B_1\to\T^2$. 
Define 
\begin{equation}\label{eq:LO}
     \Lambda_t 
     := \Lambda_t(S) 
     := \left\{[x,y+t]\,|\,[x,y]\in\p S\right\}
\end{equation}
(see Figure~\ref{fig:torus1}).  

\begin{theorem}\label{thm:2torus}
Let $S\subset\T^2$ be a closed embedded disc and 
$t\mapsto\L_t$ be the exact Lagrangian loop
defined by~(\ref{eq:LO}).  Then
$$
    \nu(\L) = \area(S).
$$
\end{theorem}

\begin{proof}
We prove that $\nu(\L)\le\area(S)$. To see this,
choose smooth functions $x,y:\R\to\R$ such that
$$
    x(\theta+1) = x(\theta),\qquad
    y(\theta+1) = y(\theta),
$$
and the map $\i_t:\R/\Z\to\T^2$ defined by 
$$
    \i_t(\theta) := [x(\theta),y(\theta)+t]
$$
is an embedding with $\i_t(\R/\Z)=\Lambda_t$. Then
$$
    \alpha_t 
    := \o(\p_t\i_t,d\i_t\cdot) 
     = -\dot xd\theta \in\Om^1(\R/\Z).
$$
Hence $\alpha_t=dh_t$ where $h_t=-x:\R/\Z\to\R$.  Hence
$$
     \|\p_t\L_t\| = \|dh_t\| = \max x-\min x
$$
and this implies
$$
     \ell(\L) = \max x-\min x.
$$
By Proposition~\ref{prop:discs} in the appendix, 
two loops $t\mapsto\Lambda_t(S)$ and $t\mapsto\Lambda_t(S')$,
associated to two embedded discs $S,S'\subset\T^2$
via~(\ref{eq:LO}), are Hamiltonian isotopic if and only if
$S$ and $S'$ have the same area.  
Now for every $\delta>0$ there
exists an embedded disc $S'$ 
(as illustrated in Figure~\ref{fig:torus2}) 
such that 
$$
     \area(S)=\area(S'),\qquad
     \max x' - \min x' < \area(S)+\delta,
$$
where $x',y':\R/\Z\to\R$ are chosen such that 
the map $\i'(\theta)=[x'(\theta),y'(\theta)]$
defines an embedding $\R/\Z\to\T^2$ whose image
is $\p S'$.  Hence the length of the loop
$t\mapsto\Lambda_t(S')$ is bounded above 
by $\area(S)+\delta$.  Thus we have proved that
$$
     \nu(\Lambda) \le \area(S).
$$

\begin{figure}[htp]
\centerline{\psfig{figure=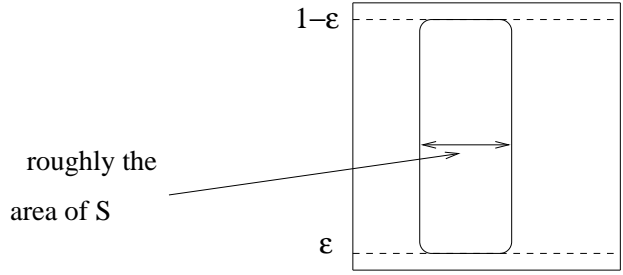}} 
\caption{{Minimizing the length}}\label{fig:torus2}
\end{figure}

\smallbreak

To show the reverse inequality let 
$t\mapsto\L'_t$ be an exact Lagrangian loop 
that is Hamiltonian isotopic to $\L$.
Then 
$$
     \L'_0=\p S',
$$
where $S'\subset\T^2$ is a smoothly 
embedded closed disc of the same area as $S$. 
Let $\psi_t:\T^2\to\T^2$ be a Hamiltonian isotopy such that
$$
     \psi_t(\L'_0)=\L'_t.
$$
We shall prove that 
\begin{equation}\label{eq:area-length}
     \area(S)\le \ell\left(\{\psi_t\}_{0\le t\le 1}\right).
\end{equation}
To see this, choose an embedded  closed discs
$\tilde S\subset\R^2$ such that 
$\pi(\tilde S)=S'$ and let 
$\tilde\psi_t:\R^2\to\R^2$ be a lift 
of $\psi_t$.  Since $\L'$ is Hamiltonian isotopic to $\L$
we have 
$
    \tilde\psi_{t+1}(\tilde S)=\tilde\psi_t(\tilde S)+(0,1)
$
and hence
$$
    \tilde\psi_1(\tilde S)\cap\tilde S = \emptyset.
$$
Let $\tilde H_t:\R^2\to\R$ be 
the Hamiltonian functions that generate $\tilde\psi_t$
and have mean value zero over the fundamental
domain $[0,1]^2$.  Choose $R>1$ such that 
$\tilde\psi_t(\tilde S)\subset B_R$ for every $t\in[0,1]$ 
and let $\beta:\R^2\to[0,1]$ be a compactly supported 
cutoff function such that $\beta|_{B_R}\equiv1$.  
Then the functions 
$$
    \hat{H}_t:=\beta\tilde H_t
$$
generate a compactly supported Hamiltonian isotopy
$\hat{\psi}_t$ of $\R^2$ that satisfies 
$$
    \hat{\psi}_1(\tilde S) \cap\tilde S = \emptyset.
$$
Now it follows from the {\it energy-capacity inequality}
in Hofer~\cite{H} that the {\it displacement energy} 
of $\tilde S$ is bounded below by the area.  
Hence
$$
     \area(\tilde S) 
     \le d(\id,\hat{\psi}_1) 
     \le \ell(\{\hat{\psi}_t\}_{0\le t\le 1})
      =  \ell(\{\psi_t\}_{0\le t\le 1}).
$$
Since 
$$
     \area(\tilde S)=\area(S')=\area(S),
$$
this proves~(\ref{eq:area-length}).  
It follows from~(\ref{eq:area-length}) 
and~(\ref{eq:length1}) that 
$$
     \area(S) \le \ell(\Lambda')
$$
for every exact Lagrangian loop $\L'$
that is Hamiltonian isotopic to $\L$.
Hence $\area(S)\le\nu(\Lambda)$.
\QED\medskip
\end{proof}

Theorem~\ref{thm:2torus} shows that the invariant
$\nu(\Lambda)$ is not necessarily invariant
under Lagrangian isotopy, but only under exact 
Lagrangian isotopy.  The techniques of proof 
are specific to the 2-dimensional case.
To establish lower bounds for our invariants
in higher dimensions we shall use existence theorems 
for pseudoholomorphic discs.


\section{Relative Gromov invariants}\label{sec:gro}

Throughout we assume that our symplectic manifold $(M,\o)$
is compact. The relative Gromov invariants of an
exact Lagrangian loop $\L\subset D\times M$
are defined in terms of holomorphic
sections of the bundle $D\times M\to D$ 
with boundary values in $\L$.
Let us denote by $\Map_\L(D,M)$ the space of smooth
functions $u:D\to M$ that satisfy
$
     u(e^{2\pi it}) \in \L_t
$
for every $t\in\R$. The {\bf Maslov class} is a function
$$
     \mu_\L:\Map_\L(D,M)\to\Z
$$
defined as follows.  Given $u\in\Map_\L(D,M)$
choose a trivialization of the tangent bundle $u^*TM$.  
Then the tangent spaces $T_{u(e^{2\pi it})}\Lambda_t$ 
define a loop of Lagrangian subspaces in $(\R^{2n},\o_0)$ 
and $\mu_\L(u)$ is defined as the 
Maslov index of this loop (cf.~\cite{RS1}).   
This integer is independent of the choice of the trivialization 
used to define it, and it depends only on the homology
class of $u$ in $H_2(D\times M,\L;\Z)$. 
We shall assume throughout that
the pair $(M,\L_0)$ is {\bf monotone}, 
i.e. there exists a $\lambda>0$ such that, 
for every smooth map $v\in\Map_{\L_0}(D,M)$,
$$
       \int_D v^*\o = \l \mu_{\L_0}(v).
$$
Here $\mu_{\L_0}$ denotes the Maslov class corresponding
to the constant loop $t\mapsto\L_0$. 
The {\bf minimal Maslov number}\label{minmas}
of the pair $(M,\L_0)$ is defined by 
$$
     N := \inf\left\{\mu_{\L_0}(v)\,|\,v:(D,\p D)\to(M,\L_0),\,
          \mu_{\L_0}(v)>0\right\}.
$$
We shall define relative Gromov invariants for
every tuple ${\bf t} = (t_1,\dots,t_k)\in\R^k$ with
$0\le t_1<\cdots<t_k<1$ and every class 
$A\in H_2(D\times M,\L;\Z)$ that satisfies
$
     n\pm\mu_\L(A) \le N-2.
$
The invariants are homology classes
$$
     \Gr^\pm_{A,{\bf t}}(\L)
     \in H_{n\pm\mu_\L(A)}(\L_{\bf t};\Z_2),
$$
where $\L_{\bf t}:=\L_{t_1}\times\cdots\times\L_{t_k}$.
These homology classes arise from
certain moduli spaces $\cM_A(\tau,\pm J)$
of (anti-)holomorphic sections of 
the bundle $D\times M$ with boundary values in $\L$
that represent the class $A$.
The points $(e^{2\pi it_1},\dots,e^{2\pi it_k})$ 
determine an evaluation map 
$$
     \ev_{\bf t}:\cM_A(\tau,\pm J)\to\L_{\bf t}
$$
and $\Gr^\pm_{A,{\bf t}}(\L)$ is defined as 
the image of the fundamental cycle of 
$\cM_A(\tau,\pm J)$ under the 
induced homomorphism on homology. 
We shall work with almost complex structures on 
$D\times M$ that are compatible with the fibration.
Every such structure is determined by a family of 
almost complex structures on $M$ and a
connection $2$-form $\tau\in\cT(\L)$.


\subsection{J-holomorphic discs}

Let $\L\subset S^1\times M$ be an exact Lagrangian loop
and $\tau\in\cT(\L)$ be a Hamiltonian connection $2$-form 
that preserves $\L$. Throughout we shall denote by
$\cJ(M,\o)$ the space of almost complex structures
on $TM$ that are compatible with $\o$.  
Let 
$
     D\to\cJ(M,\o):(x,y)\mapsto J_{x,y}
$
be a smooth family of such almost complex structures.  
Associated to the triple $(\tau,J,\L)$ there is a natural
boundary value problem for smooth functions 
$u:D\to M$:
\begin{equation}\label{eq:jhol}
     \p_xu-X_F(u) + J(\p_yu-X_G(u)) = 0,
\end{equation}
\begin{equation}\label{eq:bc}
     u(e^{2\pi it}) \in \L_t,\qquad t\in\R.
\end{equation}
Here we abbreviate $J=J_{x,y}$,
$\tau$ is given by~(\ref{eq:tau}),
$X_F=X_F(x,y,\cdot)\in\Vect(M)$ denotes
the Hamiltonian vector field of the function 
$F=F(x,y,\cdot):M\to\R$, and similarly for $X_G$.
Following Gromov~\cite{G1} we observe that the 
solutions of~(\ref{eq:jhol}) can be thought
of as pseudo-holomorphic curves in $D\times M$.

\begin{remark}\label{rmk:tildeJ}
Consider the almost complex structure
$\tilde J$ on $D\times M$ given by
$$
   \tilde J 
   = \tilde J(\tau,J)
   := \left(\begin{array}{ccc} 
     0 & -1 & 0 \\
     1 & 0 & 0 \\
     - JX_F + X_G & -X_F -JX_G & J 
     \end{array}\right).
$$
Then $u:D\to M$ is a solution of~(\ref{eq:jhol})
if and only if the function 
\begin{equation}\label{eq:tildeu}
     \tilde u(x,y)=(x,y,u(x,y))
\end{equation}
is a $\tilde{J}$-holomorphic curve in $D\times M$,
i.e.
$$
     \p_x\tilde u + \tilde{J}\p_y\tilde u=0.
$$
If $\tau$ is given by~(\ref{eq:tau}) then,
for every $\tilde\zeta=(\xi,\eta,\zeta)\in T_{x,y,z}(D\times M)$,
$$
     \tau(\tilde\zeta,\tilde J\tilde\zeta)
     = \left|\zeta-\xi X_F-\eta X_G\right|^2
       + (c-\Om_\tau)(\xi^2+\eta^2).
$$
Hence $\tilde{J}$ is tamed by $\tau$ whenever
$\tau\in\cT^+(\L)$ (see~(\ref{eq:tau-n+1})). 
If $\tau\in\cT^-(\L)$ then $\tilde J(\tau,-J)$ 
is tamed by $-\tau$. 
\end{remark}

The {\bf energy} of a solution $u$ of~(\ref{eq:jhol})
is defined by 
$$
     E(u) := \int_D\left|\p_xu-X_F(u)\right|^2\,dxdy.
$$
The next lemma shows that the solutions 
of~(\ref{eq:jhol}) and~(\ref{eq:bc}) 
that represent a given homology class 
$A\in H_2(D\times M,\L;\Z)$ satisfy a
uniform energy bound. 

\begin{lemma}\label{le:energy}
Let $u:D\to M$ be a smooth solution 
of~(\ref{eq:jhol}) and~(\ref{eq:bc})
and denote by $A\in H_2(D\times M,\L;\Z)$
the homology class represented by the map
$\tilde u:D\to D\times M$ defined by~(\ref{eq:tildeu}).
Let $c:D\to\R$ be the function in~(\ref{eq:tau}). 
Then
$$
     E(u) 
     = \inner{[\tau]}{A} 
       + \int_D \left(\Om_\tau(x,y,u)-c(x,y)\right)\, dxdy.
$$
\end{lemma}

\begin{proof}
We compute
\begin{eqnarray*}
     E(u)
&= &
     \int_D \o(\p_xu - X_F(u), 
               \p_yu - X_G(u))  \,dxdy  \\
&= &
     \int_D \biggl(
      \o(\p_xu,\p_yu) 
      - dF(u)\p_yu + dG(u)\p_xu 
      + \{F,G\}(u)\biggr)\, dxdy  \\
&= &
     \int_D \biggl(
      \o(\p_xu,\p_yu) 
      - dF(u)\p_yu + dG(u)\p_xu
     \biggr)\,dxdy  \\
&&
     +\,\int_D 
     \biggl(
     \Om_\tau(x,y,u) - (\p_yF)(u) + (\p_xG)(u)
     \biggr)\,dxdy  \\
&= &
     \int_D 
     \biggl(
     \tau(\p_x\tilde u,\p_y\tilde u) 
     - c(x,y)\biggr)\,dxdy 
     + \int_D \Om_\tau(x,y,u)\,dxdy.  
\end{eqnarray*}
This proves the lemma.
\QED\medskip
\end{proof}

Let us denote the moduli space of solutions 
of~(\ref{eq:jhol}) and~(\ref{eq:bc}) that
represent a given homology class 
$A\in H_2(D\times M,\L;\Z)$ by
$$
     \cM_A(\tau,J)
     := \left\{u:D\to M\,|\,
        u\mbox{ satisfies }
        (\ref{eq:jhol})\mbox{ and }(\ref{eq:bc}),\,
        [\tilde u]=A\right\}.
$$
We shall prove that, for a generic pair $(\tau,J)$,
this space is a smooth manifold of dimension
$
     n + \mu_\L(A).
$
Moreover, if the pair $(M,\L_0)$ is monotone with
minimal Maslov number $N$ and $n+\mu_\L(A)<N$,
we shall prove that $\cM_A(\tau,J)$ is compact, 
again for a generic pair $(\tau,J)$. 
The key tool for establishing compactness is
the energy bound of Lemma~\ref{le:energy}.
Under these asumptions the moduli spaces can be
used to define Gromov invariants of $\L$. 
The significance of these invariants for 
exact Lagrangian loops lies in the following 
observation.

\begin{lemma}\label{le:eps+}
Let $\L$ be an exact Lagrangian 
loop and $A\in H_2(D\times M,\L;\Z)$. 
Suppose that for every $\tau\in\cT^+(\L)$ 
there exists a $J$ such that
$
     \cM_A(\tau,J) \ne\emptyset.
$
Then 
$$
     \eps^+(\tau_0,\L) + \inner{[\tau_0]}{A} \ge 0 
$$
for every $\tau_0\in\cT(\L)$.
\end{lemma}

\begin{proof}
Let $\tau\in\cT^+(\L)$ and $u\in\cM_A(\tau,J)$. 
Let $\tilde u:D\to D\times M$ be given 
by~(\ref{eq:tildeu}). Then $\tilde u$ is a
$\tilde J(\tau,J)$-holomorphic curve.
By Remark~\ref{rmk:tildeJ}, 
$\tilde{J}(\tau,J)$ is tamed by $\tau$.  
Hence
$$
    0 < \int_D {\tilde u}^*\tau
      = \inner{[\tau]}{A}
      = \inner{[\tau_0]}{A} + s(\tau,\tau_0).
$$
The infimum of the numbers on the right
is $\inner{[\tau_0]}{A}+\eps^+(\tau_0,\L)$.
This proves the lemma.
\QED\medskip
\end{proof}

A similar estimate for $\eps^-(\tau_0,\L)$ can be 
obtained by studying anti-ho\-lo\-mor\-phic curves.
These are solutions of the equation
\begin{equation}\label{eq:jhol-}
     \p_xu-X_F(u) - J(\p_yu-X_G(u)) = 0,
\end{equation}
that satisfy the same boundary condition~(\ref{eq:bc}).
Let us denote the moduli space of solutions by 
$\cM_A(\tau,-J)$

\begin{lemma}\label{le:eps-}
Let $\L$ be an exact Lagrangian 
loop and $A\in H_2(D\times M,\L;\Z)$. 
Suppose that for every $\tau\in\cT^-(\L)$ 
there exists a $J$ such that
$
     \cM_A(\tau,-J) \ne\emptyset.
$
Then 
$$
     \eps^-(\tau_0,\L) + \inner{[\tau_0]}{A} \le 0 
$$
for every $\tau_0\in\cT(\L)$.
\end{lemma}

\begin{proof}
Let $\tau\in\cT^-(\L)$ and $u\in\cM_A(\tau,-J)$. 
Let $\tilde u:D\to D\times M$ be given 
by~(\ref{eq:tildeu}). Then $\tilde u$ is a
$\tilde J(\tau,-J)$-holomorphic curve.
By Remark~\ref{rmk:tildeJ}, 
$\tilde{J}(\tau,-J)$ is tamed by $-\tau$.  
Hence
$$
    0 > \int_D {\tilde u}^*\tau
      = \inner{[\tau]}{A}
      = \inner{[\tau_0]}{A} + s(\tau,\tau_0).
$$
The supremum of the numbers on the right
is $\inner{[\tau_0]}{A}+\eps^-(\tau_0,\L)$.
This proves the lemma.
\QED\medskip
\end{proof}


\subsection{Fredholm theory}

In this subsection we examine the moduli spaces
$\cM_A^\pm(\tau,J)$ in more detail and show that,
for a generic $J$, these spaces are smooth manifolds
of the predicted dimensions $n\pm\mu_\L(A)$. 
The arguments are standard (cf.~\cite{FHS,MS2})
and we shall only outline the main points.
Fix an exact Lagrangian loop $\L\subset S^1\times M$,
a homology class $A\in H_2(D\times M,\L;\Z)$,
and a constant $p>2$.  Consider the Banach manifold 
$$
     \cB = W^{1,p}_{\L,A}(D,M)
$$
of all functions $u:D\to M$ of class $W^{1,p}$ that satisfy the 
boundary condition~(\ref{eq:bc}) and represent the class $A$. 
There is a natural vector bundle 
$
     \cE\to\cB
$
with fibres 
$$
     \cE_u=L^p(D,u^*TM)
$$
and the left hand sides
of~(\ref{eq:jhol}) and~(\ref{eq:jhol-}) define 
Fredholm sections 
$
     \cF^\pm:\cB\to\cE
$
given by
$$
     \cF^\pm(u) 
     := \cF(u;\tau,\pm J) 
     := \p_xu-X_F(u) \pm J(\p_yu-X_G(u)).
$$
The moduli spaces $\cM_A(\tau,\pm J)$ 
are the zero sets of these sections.
The tangent space 
$$
     T_u\cB = W^{1,p}_\L(D,u^*TM)
$$
consists of all vector fields 
$\xi\in W^{1,p}(D,u^*TM)$ along $u$ which
are of class $W^{1,p}$ and satisfy the boundary
condition
$
     \xi(e^{2\pi it}) \in T_{u(e^{2\pi it})}\L_t.
$
The vertical differential of $\cF^\pm$
at a zero $u\in\cM^\pm_A(\tau,J)$ is the linear operator
$$
     D_u^\pm
     = D\cF^\pm(u):W^{1,p}_\L(D,u^*TM)
     \to L^p(D,u^*TM)
$$
given by
\begin{equation}\label{eq:Du}
     D_u^\pm\xi
     = \Nabla{x}\xi - \Nabla{\xi}X_F(u) 
       \pm J(\Nabla{y}\xi - \Nabla{\xi}X_G(u))  
       \pm (\Nabla{\xi}J)(\p_yu-X_G(u)).
\end{equation}
Here $\nabla$ denotes 
the Levi-Civita connection of the Riemannian metric
$$
     \inner{\cdot}{\cdot}_{x,y}=\o(\cdot,J_{x,y}\cdot)
$$
and thus depends on $x+iy\in D$. 
The expression $\nabla X_F$ denotes the covariant 
derivative of $X_F=X_{F_{x,y}}$ with respect 
to the Levi-Civita connection at the point $x+iy$. 
The next theorem follows from the Riemann-Roch theorem
for discs (see for example~\cite{RS2} 
for a recent exposition) and the infinite dimensional 
implicit function theorem (see for example~\cite[Appendix~B]{S}).
The proof is standard (see for example~\cite{MS2})
and will be omitted.

\begin{theorem}\label{thm:fred}
For every $u\in W^{1,p}_{\L,A}(D,M)$ the operators
$D_u^\pm$ defined by~(\ref{eq:Du}) are Fredholm
and their indices are
$$
     {\rm index} D_u^\pm = n\pm\mu_\L(u).
$$
If $D_u^\pm$ is surjective for every 
$u\in\cM_A(\tau,\pm J)$ then $\cM_A(\tau,\pm J)$
is a smooth manifold of dimension 
$$
     \dim\cM_A(\tau,\pm J) = n\pm\mu_\L(A).
$$
\end{theorem}

Fix an exact Lagrangian loop $\L$ and
a connection $2$-form $\tau\in\cT(\L)$. 
Let us denote by $\cJ(D;M,\o)$ the space
of all smooth families of almost complex 
structures $J:D\to\cJ(M,\o)$. 
Such a family $J\in\cJ(D;M,\o)$ is called {\bf regular}
for~(\ref{eq:jhol}) if $D_u^+$ is surjective 
for every $A\in H_2(D\times M,\L;\Z)$ and every
$u\in\cM_A(\tau,J)$. Similarly, 
$J\in\cJ(D;M,\o)$ is called {\bf regular}
for~(\ref{eq:jhol-}) if $D_u^-$ is surjective 
for every $A\in H_2(D\times M,\L;\Z)$ and every
$u\in\cM_A(\tau,-J)$.  We shall denote 
set of all families of almost complex structures 
that are regular for~(\ref{eq:jhol}), 
respectively~(\ref{eq:jhol-}), by
$$
     \Jreg^\pm(\tau,\L)\subset \cJ(D;M,\o).
$$
The proof of the next theorem is a standard application of the 
Sard-Smale theorem (cf.~\cite{MS2}) and will be omitted. 

\begin{theorem}\label{thm:sard}
The sets $\Jreg^\pm(\tau,\L)$ are of the second category 
in $\cJ(D;M,\o)$ in the sense of Baire, 
i.e. they are countable intersections
of open and dense subsets of $\cJ(D;M,\o)$.
In particular, they are dense.
\end{theorem}

Let $\tau_0,\tau_1\in\cT(\L)$
and choose regular families of almost complex structures
$$
     J_0\in\Jreg^\pm(\tau_0,\L),\qquad
     J_1\in\Jreg^\pm(\tau_1,\L).
$$
By Theorem~\ref{thm:fred},
the spaces $\cM_A(\tau_0,\pm J_0)$ and 
$\cM_A(\tau_1,\pm J_1)$ are smooth manifolds 
of the same dimension. These manifolds are cobordant. 
To construct a cobordism choose a smooth path
$
     [0,1]\to\cT(\L):\l\mapsto\tau_\l
$
that connects $\tau_0$ to $\tau_1$.  Let us denote by 
$$
     \cJ = \cJ([0,1]\times D,J_0,J_1;M,\o)
$$
the space of smooth homotopies 
$
     [0,1]\to\cJ(D;M,\o):\l\mapsto J_\l
$
that connect $J_0$ to $J_1$. 
Given $\{J_\l\}\in\cJ$ denote
$$
     \cW_A(\{\tau_\l\},\{\pm J_\l\})
     = \left\{(\l,u)\,|\,0\le\l\le1,\,
       u\in\cM_A(\tau_\l,\pm J_\l)\right\}.
$$
A homotopy $\{J_\l\}\in\cJ$ is called {\bf regular} if,
for every $A\in H_2(D\times M,\L;\Z)$ and every 
pair $(\l,u)\in\cW_A(\{\tau_\l\},\{\pm J_\l\})$,
$$
     \im D^\pm_{\l,u} + \R \xi^\pm_{\l,u} 
     = L^p(D,u^*TM).
$$
Here $D^\pm_{\l,u}$ is defined by~(\ref{eq:Du})
with $\tau$ and $J$ replaced by $\tau_\l$ and $J_\l$,
respectively, and $\xi^\pm_{\l,u}\in L^p(D,u^*TM)$
is given by
$$
     \xi^\pm_{\l,u} 
     := X_{\p_\l F_\l}(u)
        \pm J_\l(u)X_{\p_\l G_\l}(u)
        \mp \p_\l J_\l(u)(\p_yu-X_{G_\l}(u)).
$$
The set of all regular homotopies will be denoted by
$$
     \Jreg^\pm(\{\tau_\l\},J_0,J_1,\L)\subset\cJ.
$$
The proof of the next theorem
is again standard and will be omitted.

\begin{theorem}\label{thm:homotopy}
Let $[0,1]\to\cT(\L):\l\mapsto\tau_\l$
be a smooth family of connection $2$-forms
and suppose that 
$
     J_0\in\Jreg^\pm(\tau_0,\L)
$
and
$
     J_1\in\Jreg^\pm(\tau_1,\L).
$
Then the sets $\Jreg^\pm(\{\tau_\l\},J_0,J_1,\L)\subset\cJ$ 
are of the second category in the sense of Baire.
Moreover, if $\{J_\l\}\in\Jreg^\pm(\{\tau_\l\},J_0,J_1,\L)$ 
then $\cW_A(\{\tau_\l\},\{\pm J_\l\})$
is a smooth manifold of dimension
$$
     \dim\cW_A(\{\tau_\l\},\{\pm J_\l\})
     = n\pm\mu_\L(A)+1
$$ 
with boundary
$$
     \p\cW_A(\{\tau_\l\},\{\pm J_\l\})
     = \cM_A(\tau_0,\pm J_0)\cup \cM_A(\tau_1,\pm J_1).
$$
\end{theorem}


\subsection{Compactness}

\begin{theorem}\label{thm:compact}
Let $\L\subset S^1\times M$ be an exact Lagrangian
loop and suppose that the pair $(M,\L_0)$ is monotone.
Let $A\in H_2(D\times M,\L;\Z)$ and denote by
$N\in\N$ the minimal Maslov number 
of the pair $(M,\L_0)$.

\smallskip
\noindent{\bf (i)}
If 
$$
    n\pm\mu_\L(A)\le N-1
$$
then $\cM_A(\tau,\pm J)$ is compact for every 
$\tau\in\cT(\L)$ and every $J\in\Jreg^\pm(\tau,\L)$.

\smallskip
\noindent{\bf (ii)}
If 
$$
    n\pm\mu_\L(A)\le N-2
$$
then $\cW_A(\{\tau_\l\},\{\pm J_\l\})$ is compact for 
every smooth path $[0,1]\to\cT(\L):\l\mapsto\tau_\l$,
every $J_0\in\Jreg^\pm(\tau_0,\L)$,
every $J_1\in\Jreg^\pm(\tau_1,\L)$,
and every regular homotopy 
$\{J_\l\}\in\Jreg^\pm(\{\tau_\l\},J_0,J_1,\L)$. 
\end{theorem}

The proof of Theorem~\ref{thm:compact} relies
on the following theorem about Gromov compactness
for $J$-holomorphic discs.  This result is implicitly
contained in Gromov's original paper~\cite{G1}
and has been folk knowledge since then.
However, the full details of the proof have not so
far appeared in the literature.  They were recently 
carried out by Frauenfelder~\cite{F}
in his Diploma thesis.  In his thesis Frauenfelder also 
discusses the corresponding notion of {\it stable maps}
for pseudoholomorphic discs. 

\begin{theorem}[Gromov]\label{thm:gromov}
Let $(\tau^\nu,J^\nu)\in\cT(\L)\times\cJ(D;M,\o)$
be a sequence that converges in the $\Cinf$-topology
to $(\tau,J)\in\cT(\L)\times\cJ(D;M,\o)$.
Let $A\in H_2(D\times M,\L;\Z)$ and 
$u^\nu\in\cM_A(\tau^\nu,\pm J^\nu)$.
If $u^\nu$ has no $\Cinf$-convergent subsequence
then there exist 
\begin{description}
\item[(i)]
finitely many points $(x_i,y_i)\in D$ and
maps $v_i:S^2\to M$, $i=1,\dots,k$,
\item[(ii)]
finitely many points $t_j\in\R$ and 
maps $w_j:D\to M$,
$j=1\dots,\ell$,
\item[(iii)]
a map $u_0:D\to M$,
\end{description}
such that $v_i$ is a nonconstant 
$J_{x_i,y_i}$-(anti)holomorphic 
sphere for $i=1,\dots,k$, 
$w_j$ is a nonconstant 
$J_{e^{2\pi it_j}}$-(anti)holomorphic disc 
with $w_j(\p D)\subset\L_{t_j}$ for $j=1,\dots,\ell$, 
$u_0\in\cM_{A_0}(\tau,\pm J)$ for some 
$A_0\in H_2(D\times M,\L;\Z)$, and
\begin{equation}\label{eq:A}
     A = A_0 + \sum_{i=1}^k[v_i] + \sum_{j=1}^\ell[w_j].
\end{equation}
Here $[v_i]$ and $[w_j]$ denote the induced homology
classes in $H_2(D\times M,\L;\Z)$ and one of the 
integers $k$ and $\ell$ is nonzero.
\end{theorem}

\begin{remark}\label{rmk:gromov}
{\bf (i)}
Let $\tilde M$ be a compact manifold
and $\tilde L\subset\tilde M$ be a compact 
submanifold of half the dimension.
Suppose that $\tilde\o^\nu$ is a sequence of symplectic 
forms on $M$ that converges to $\tilde\o$ in the $\Cinf$-topology
such that $\tilde L$ is a Lagrangian submanifold
of $(\tilde M,\tilde\o^\nu)$ for every $\nu$.  
Suppose that $\tilde J^\nu$ is a sequence 
of $\tilde\o^\nu$-tame almost complex structures 
on $\tilde M$ that converges
in the $\Cinf$-topology to $\tilde J$.  
In~\cite{F} Frauenfelder proves, in particular, 
that a sequence of $\tilde J^\nu$-holomorphic discs
$\tilde u^\nu:(D,\p D)\to(\tilde M,\tilde L)$
that represent a fixed homology class 
$A\in H_2(\tilde M,\tilde L;\Z)$
has a subsequence that converges (in a precisely defined sense)
to a tree consisting of $\tilde J$-holomorphic spheres in $\tilde M$ 
and $\tilde J$-holomorphic discs in $\tilde M$ with boundary in $\tilde L$
such that the sum of their homology classes in
$H_2(\tilde M,\tilde L;\Z)$ is equal to $A$.
The techniques in~\cite{F} are an adaptation of those
in Hofer--Salamon~\cite{HS} for holomorphic spheres 
to the case of holomorphic discs.

\smallskip
\noindent{\bf (ii)}
The moduli space $\cM_A(\tau,\pm J)$ does not
depend on the function $c:D\to M$ in~(\ref{eq:tau}).
Hence we may assume without loss of generality that the 
connection forms $\tau^\nu$ in Theorem~\ref{thm:gromov}
lie in $\cT^\pm(\L)$. 
Under this assumption the manifold $\tilde M=D\times M$, 
the submanifold $\tilde L=\L$,
the symplectic forms $\tilde\o^\nu=\pm\tau^\nu$,
the almost complex structures 
$\tilde J^\nu=\tilde J(\tau^\nu,\pm J^\nu)$
defined in Remark~\ref{rmk:tildeJ},
and the functions $\tilde u^\nu$ given by~(\ref{eq:tildeu})
satisfy the requirements of~(i). 

\smallskip
\noindent{\bf (iii)}
Theorem~\ref{thm:gromov} follows from~(i) and~(ii)
since each bubble in the limit curve is contained 
in a fibre of the (trivial) fibration $D\times M$. 
To see this, note that each curve $v_i$ appears as
the limit of a sequence 
$$
       v^\nu_i(x,y)
       = u^\nu(x^\nu_i+\eps^\nu x,y^\nu_i+\eps^\nu y),
$$
where $x^\nu_i\to x_i$, $y^\nu_i\to y_i$, $\eps^\nu\to 0$,
and 
$$
      \lim_{\nu\to\infty}
      \frac{\eps^\nu}{1-\sqrt{(x^\nu_i)^2+(y^\nu_i)^2}}
      = 0.
$$ 
One can show that, after passing to a suitable subsequence, 
the sequence $v^\nu_i$ converges to $v_i$ 
in the $\Cinf$-topology on the complement of some
finite set. The functions $v^\nu_i$ satisfy
$$
      \p_xv^\nu_i - \eps^\nu X_{F^\nu}
      + J^\nu(\p_yv^\nu_i - \eps^\nu X_{G^\nu})
      = 0,
$$
where $X_{F^\nu}$, $X_{G^\nu}$, and $J^\nu$ 
are evaluated at the point
$(x^\nu_i+\eps^\nu x,y^\nu_i+\eps^\nu y,v^\nu_i)$.
It follows that the limit curve $v_i$
extends to a $J_{x_i,y_i}$-holomorphic sphere. 
The holomorphic discs $w_j$ appear as similar limits with
$x_j+iy_j=e^{2\pi it_j}$ and 
$$
      \lim_{\nu\to\infty}
      \frac{\eps^\nu}{1-\sqrt{(x^\nu_j)^2+(y^\nu_j)^2}}
      > 0.
$$ 
A similar argument as above, with coordinates on the 
upper halfplane, then shows that the limit curve $w_j$ 
is a $J_{e^{2\pi it_j}}$-holomorphic disc with boundary
values in $\L_{t_j}$.  

\smallskip
\noindent{\bf (iv)}
The limit curve in~(i) is a 
stable map consisting of $\tilde J$-holomorphic discs
and spheres.  For closed curves this concept is 
due to Kontsevich~\cite{K}.  Some of the components
of the stable map may be constant.  However,
these do not contribute to the homology class
and can be neglected for our purposes. 
If the original sequence $\tilde u^\nu$ does not have 
a $\Cinf$-convergent subsequence, then the limit 
curve has more than one nonconstant component.
This shows that in Theorem~\ref{thm:gromov} either
$k$ or $\ell$ is nonzero.
\end{remark}

\medskip
\noindent{\bf Proof of Theorem~\ref{thm:compact}:}
We prove~(ii).  Suppose, by contradiction, that 
$\cW_A(\{\tau_\l\},\{\pm J_\l\})$ is not compact.  
Then there exists a sequence 
$$
     (\l^\nu,u^\nu)\in\cW_A(\{\tau_\l\},\{\pm J_\l\})
$$
that has no convergent subsequence. 
We may assume without loss
of generality that $\l^\nu$ converges to $\l_0$. 
Then, by Theorem~\ref{thm:gromov},
there exist nonconstant $J_{\l_0;x_i,y_i}$-(anti)holomorphic 
spheres $v_i:S^2\to M$ for $i=1,\dots,k$, 
nonconstant $J_{\l_0;e^{2\pi it_j}}$-(anti)holomorphic discs
$w_j:(D,\p D)\to(M,L_{t_j})$ for $j=1,\dots,\ell$, 
and an element $u_0\in\cM_{A_0}(\tau_{\l_0},\pm J_{\l_0})$
for some $A_0\in H_2(D\times M,\L;\Z)$ such that~(\ref{eq:A})
is satisfied.  Since the pair $(M,\L_t)$ is 
monotone with minimal Maslov number $N$ for every
$t$ we have 
$$
      \pm\mu_\L(v_i) \ge N,\qquad
      \pm\mu_\L(w_j) \ge N
$$
for $i=1,\dots,k$ and $j=1,\dots,\ell$.
Since either $k$ or $\ell$ is nonzero this implies
\begin{eqnarray*}
      n\pm\mu_\L(A) 
&= &
      n\pm\mu_\L(A_0) 
        \pm \sum_{i=1}^k\mu_\L(v_i)
        \pm \sum_{j=1}^\ell\mu_\L(w_j) \\
&\ge &
      n\pm\mu_\L(A_0)+N. 
\end{eqnarray*}
Since $\{J_\l\}\in\Jreg^\pm(\{\tau_\l\},J_0,J_1,\L)$,
the moduli space $\cW_{A_0}(\{\tau_\l\},\{\pm J_\l\})$ 
is a smooth manifolds of dimension
\begin{eqnarray*}
      \dim\,\cW_{A_0}(\{\tau_\l\},\{\pm J_\l\})
&= &
      n\pm\mu_\L(A_0)+1  \\
&\le &
      n\pm\mu_\L(A)+1-N  \\
&< &
      0.
\end{eqnarray*}
Hence 
$$
     \cW_{A_0}(\{\tau_\l\},\{\pm J_\l\})=\emptyset,
$$
in contradiction to the fact that 
$$
     (\l_0,u_0)\in\cW_{A_0}(\{\tau_\l\},\{\pm J_\l\}).
$$
Thus we have proved~(ii).  The proof of~(i) 
is almost word by word the same and will be 
left to the reader. 
\QED


\subsection{Gromov invariants}

Fix an exact Lagrangian loop $\L\subset S^1\times M$
and a class $A\in H_2(D\times M,\L;\Z)$.
Throughout we shall assume that the pair $(M,\L_0)$
is monotone and 
\begin{equation}\label{eq:nAN}
        n\pm\mu_\L(A)\le N-2,
\end{equation}
where $N\in\N$ denotes the minimal Maslov number of 
the pair $(M,\L_0)$. 
Fix a tuple ${\bf t}=(t_1,\dots,t_k)\in\R^k$ 
such that $0\le t_1<\cdots<t_k<1$ and denote
$$
    \L_{\bf t}:=\L_{t_1}\times\cdots\times\L_{t_k}.
$$
For $\tau\in\cT(\L)$ and $J\in\cJ(D;M,\o)$
we define
$
     \ev_{\bf t}:\cM_A(\tau,\pm J)\to\L_{\bf t}
$
by 
$$
     \ev_{\bf t}(u)
     := (u(e^{2\pi it_1}),\dots,u(e^{2\pi it_k})).
$$
If $J\in\Jreg^\pm(\tau,\L)$ then, 
by Theorems~\ref{thm:fred} and~\ref{thm:compact},
the moduli space $\cM^\pm_A(\tau,J)$ is a compact 
smooth manifolds(without boundary) 
of dimension $n\pm\mu_\L(A)$.
It is not necessarily orientable.
Let 
$$
     [\cM_A(\tau,\pm J)]
     \in H_{n\pm\mu_\L(A)}(\cM_A(\tau,\pm J);\Z_2)
$$
denote the fundamental cycle.
The Gromov invariants
are defined by 
\begin{equation}\label{eq:gromov}
     \Gr^\pm_{A,{\bf t}}(\L)
     := {\ev_{\bf t}}_*[\cM^\pm_A(\tau,J)]
     \in H_{n\pm\mu_\L(A)}(\L_{\bf t};\Z_2).
\end{equation}

\begin{lemma}\label{le:gromov}
The homology classes 
$
     \Gr^\pm_{A,{\bf t}}(\L)
     \in H_{n\pm\mu_\L(A)}(\L_{\bf t};\Z_2)
$
are independent of the choices of the connection $2$-form
$\tau\in\cT(\L)$ and the almost complex structure 
$J\in\Jreg^\pm(\tau,\L)$ used to define them. 
\end{lemma}

\begin{proof}
Theorems~\ref{thm:homotopy} and~\ref{thm:compact}~(ii).
\QED
\end{proof}

\begin{corollary}\label{cor:gromov}
Let $A^\pm\in H_2(D\times M,\L;\Z)$ satisfy~(\ref{eq:nAN})
and suppose that 
$$
     \Gr^\pm_{A^\pm,{\bf t}^\pm}(\L)\ne 0
$$
for some ${\bf t}^\pm$.  Then 
$$
     \eps^+(\tau,\L) \ge -\inner{[\tau]}{A^+},\qquad
     \eps^-(\tau,\L) \le -\inner{[\tau]}{A^-}
$$
for every $\tau\in\cT(\L)$. 
\end{corollary}

\begin{proof}
Lemmata~\ref{le:eps+} and~\ref{le:eps-}.
\QED
\end{proof}


\section{Complex projective space}\label{sec:cx}
 
In this section we shall use the Gromov invariants
to compute the K-area of certain exact Lagrangian 
loops in $\C P^n$. The archetypal example is the half turn
of a great circle in the 2-sphere. An explicit computation
shows that the Hofer length of this loop is $1/2$.
We shall use Corollary~\ref{cor:gromov} and
Theorems~\ref{thm:nu-chi} and~\ref{thm:eps-nu} to show
that this loop minimizes the Hofer length in its
Hamiltonian isotopy class. 


\subsection{Rotations of real projective space}

Consider the complex projective space 
$$
     M=\C P^n
$$
equipped with symplectic form $\o$ 
that is induced by the Fubini-Study 
metric and satisfies the normalization
condition
$$
     \int_{\C P^n}\o^n=1.
$$ 
Let $L=\R P^n$ and fix an integer $k\in\{1,\dots,n\}$.
As in the introduction, we consider the exact Lagrangian loop
\begin{equation}\label{eq:tLk}
     \L := \bigcup_{t\in\R}
     \{e^{2\pi it}\}\times\psi_t(\R P^n),
\end{equation}
where 
$$
     \psi_t([z_0:\cdots:z_n])
     = ([z_0:e^{\pi it}z_1:\cdots:e^{\pi it}z_k:z_{k+1}:\cdots:z_n]).
$$
The Hamiltonian isotopy $\psi_t$
is generated, via~(\ref{eq:ham}), by the 
the time independent Hamiltonian function
$H_t=H:\C P^n\to\R$ given by
\begin{equation}\label{eq:H}
     H([z_0:...:z_n]) 
     = \frac{k}{2n+2}
       - \frac{|z_1|^2+\cdots+|z_k|^2}{2(|z_0|^2+\cdots+|z_n|^2)}.
\end{equation}
This function has mean value zero and Hofer norm
$$
     \|H\|
     = \max H - \min H
     = \frac12.
$$
Since $H$ attains its maximum and its minimum
on $\L_t=\psi_t(\R P^n)$ it follows that 
$
    \ell(\L) = 1/2.
$


\subsection{The Maslov index}

We prove that the minimal Maslov number 
of the pair $(\C P^n,\R P^n)$ is
\begin{equation}\label{eq:N}
     N = n+1.
\end{equation}
For $n=1$ this is obvious. For $n>1$ consider the 
homology exact sequence of the pair $(\C P^n,\R P^n)$.
It has the form
$$
     0\to H_2(\C P^n;\Z)\to H_2(\C P^n,\R P^n;\Z)
      \to H_1(\R P^n;\Z)\to 0.
$$
Now $\R P^n$ decomposes the line $\C P^1\subset\C P^n$ 
into two discs that represent the same homotopy class 
in $\pi_2(\C P^n,\R P^n)$. Hence there is an element
$A\in H_2(\C P^n,\R P^n;\Z)$ such that $2A$ is equal
to the image of the generator under the homomorphism
$$
     \Z\cong H_2(\C P^n;\Z)\to H_2(\C P^n,\R P^n;\Z).
$$
This implies that $A$ is the generator of 
$H_2(\C P^n,\R P^n;\Z)\cong\Z$. Since the Maslov class
of $2A\in\pi_2(\C P^n,\R P^n)$ is equal to 
$2\inner{c_1(T\C P^n)}{[\C P^1]}=2n+2$
we have proved~(\ref{eq:N}). 

\begin{lemma}\label{le:maslov}
Let $(M,\o)$ be a symplectic manifold 
and $\L\subset S^1\times M$ be an 
exact Lagrangian loop such that $(M,\L_0)$
is a monotone pair with minimal Maslov
number $N$. Then 
$$
     \mu_\L(u_1)\equiv \mu_\L(u_0)\mbox{ {\rm mod} }N
$$
for all $u_0,u_1\in\Map_\L(D,M)$.
\end{lemma}

\begin{proof}
If $u_0(e^{2\pi it})=u_1(e^{2\pi it})$ for every $t\in\R$
then $u_0$ (with reversed orientation)
and $u_1$ form a sphere and 
the difference $\mu_\L(u_1)-\mu_\L(u_0)$ is equal
twice the first Chern number of this sphere.  
Hence the difference of the Maslov numbers 
is an even multiple of $N$.  This continues to hold 
whenever $u_0|_{\p D}$ is homotopic to $u_1|_{\p D}$
as a section of the bundle $\L\to S^1$. 
For any two maps $u_0,u_1\in\Map_\L(D,M)$ there exists
a smooth function $v:(D,\p D)\to(M,\L_0)$ such that
$v(-1)=u_0(1)$ and the connected sum $u_0\#v$ 
is homotopic to $u_1$ along the boundary.  
Hence, by what we have just proved,
$$
     \mu_\L(u_1)-\mu_\L(u_0)-\mu_{\L_0}(v)
     \in 2N\Z.
$$
Since $\mu_{\L_0}(v)$ is an integer 
multiple of $N$, the lemma is proved.
\QED\medskip
\end{proof} 

Returning to the loop $\L\subset S^1\times\C P^n$
we observe that~(\ref{eq:H}) 
is a Morse-Bott function with critical manifolds
$$
     C^+ := \left\{[0:z_1:\cdots:z_k:0:\cdots:0]\,|\,
           (z_1,\dots,z_k)\in\C^k\setminus\{0\}\right\},
$$
$$
     C^- := \left\{[z_0:0:\cdots:0:z_{k+1}:\cdots:z_n]\,|\,
           (z_0,z_{k+1},\dots,z_n)\in\C^{n-k+1}\setminus\{0\}\right\}.
$$
Note that $H$ attains its minimum $(k-n-1)/(2n+2)$ 
on $C^+$ and its maximum $k/(2n+2)$ on $C^-$.
Moreover, $C^\pm\cap\R P^n\subset\L_t$ for every $t$.
Let us denote by 
$$
     A^\pm\in H_2(D\times\C P^n,\L;\Z)
$$
the homology classes represented by the constant 
functions $D\to\C P^n$ with values in $C^\pm\cap\R P^n$. 
The next lemma shows that $\L$ 
has Maslov index $k\in\Z_{n+1}$
as claimed in the introduction (see~(\ref{eq:mu-Lk})).  
It also shows that the homology classes 
$A^\pm\in H_2(D\times\C P^n,\L;\Z)$
satisfy the condition~(\ref{eq:nAN})
for the definition of the Gromov invariants.

\begin{lemma}\label{le:mu}
\begin{equation}\label{eq:mu}
     \mu_\L(A^+) = k-1-n,\qquad
     \mu_\L(A^-) = k.
\end{equation}
\end{lemma}

\begin{proof}
In the case of $A^-$, 
consider the constant function 
$
     u(x,y)\equiv p := [1:0:\cdots:0]. 
$
Then a trivialization of the pullback tangent bundle 
$u^*T\C P^n$ is determined by the coordinate chart
$
     [z_0:\cdots:z_n]
     \mapsto(z_1/z_0,\dots,z_n/z_0).
$
In these coordinates the Hamiltonian flow is
$
     \zeta\mapsto
     (e^{\pi i t}\zeta_1,\dots,e^{\pi i t}\zeta_k,
     \zeta_{k+1},\dots,\zeta_n).
$
Since 
$
     T_p\L_0\cong\R^n\subset\C^n\cong T_p\C P^n,
$
we see that the Maslov index of the loop
$t\mapsto T_p\L_t$ is equal to $k$.
This proves the second equation in~(\ref{eq:mu})
and the first follows from a similar argument.
\QED
\end{proof}


\subsection{Computation of the Gromov invariants}

Since $N=n+1$ it follows from Lemma~\ref{le:mu}
that the classes $A^\pm$ satisfy~(\ref{eq:nAN}) and hence
the requirements of Theorem~\ref{thm:compact}.
The next theorem shows that the Gromov invariants
$\Gr_{A^\pm,0}^\pm(\L)$ are nonzero. 
Here the subscript $0$ corresponds to the choice
${\bf t}=t_1=0$ for the evaluation map.

\begin{theorem}\label{thm:cpn}
$$
     \Gr_{A^+,0}^+(\L) = [\R P^{k-1}]\in H_{k-1}(\R P^n;\Z_2),
$$
$$
     \Gr_{A^-,0}^-(\L) = [\R P^{n-k}]\in H_{n-k}(\R P^n;\Z_2).
$$
\end{theorem}

\begin{proof}
Let $\tau\in\cT(\L)$ be given by~(\ref{eq:tau}) 
with $c=0$ and 
$$
     F_{x,y} = \frac{-\sin(2\pi t)\rho(r)}{2\pi r}H,\qquad
     G_{x,y} = \frac{\cos(2\pi t)\rho(r)}{2\pi r}H,
$$
where $re^{2\pi it}=x+iy$ and $H$ is given by~(\ref{eq:H}). 
As in~(\ref{eq:FG}), $\rho:[0,1]\to[0,1]$ is a smooth 
nondecreasing cutoff function such that $\rho(r)=0$ 
for $r$ near $0$ and $\rho(r)=1$ for $r$ near $1$.
The formula 
\begin{equation}\label{eq:Omtau}
     \Om_\tau(re^{2\pi it},z) = -\frac{\dot\rho(r)}{2\pi r}H(z)
\end{equation}
for $z\in\C P^n$ shows that $\Om_\tau(x,y,z)\ge 0$ for $z\in C^+$
and $\Om_\tau(x,y,z)\le0$ for $z\in C^-$. 
By~(\ref{eq:Omtau}) and Lemma~\ref{le:energy} 
with $c=0$ and $E(u)=0$,
\begin{equation}\label{eq:tau-A}
     \inner{[\tau]}{A^+} = \frac{k-1-n}{2n+2},\qquad
     \inner{[\tau]}{A^-} = \frac{k}{2n+2}.
\end{equation}
The explicit formulae for $F$ and $G$ show that 
$C^\pm$ consist entirely of critical 
points of $F_{x,y}$ and $G_{x,y}$ for all $x+iy\in D$. 
This shows that the constant functions $u:D\to\C P^n$
with values in $C^+\cup C^-$ are horizontal for the 
symplectic connection determined by $\tau$. In explicit terms
$\p_xu=X_F(u)$ and $\p_yu=X_G(u)$. Hence these constant 
functions satisfy both equations~(\ref{eq:jhol}) 
and~(\ref{eq:jhol-}) for every $J\in\cJ(D;\C P^n,\o)$.
The constant functions with values in 
$(C^+\cup C^-)\cap\R P^n$ satisfy 
in addition the boundary condition~(\ref{eq:bc}). 
The formula~(\ref{eq:mu}) shows that
the constant solutions with values in $C^+\cap\R P^n$ 
and those with values in $C^-\cap\R P^n$ 
represent different homology classes. 

We prove that, for every $J\in\cJ(D;\C P^n,\o)$,
\begin{equation}\label{eq:constant}
     \cM_{A^+}(\tau,J)
     = \left\{u:D\to C^+\cap\R P^n\,|\,du=0\right\}.
\end{equation}
To see this, let $u\in\cM_{A^+}(\tau,J)$.  
Then, by Lemma~\ref{le:energy} and~(\ref{eq:Omtau}),
\begin{eqnarray*}
     0
&\le &
     E(u)  \\
&= &
     \inner{[\tau]}{A^+}
     + \int_D \Om_\tau(x,y,u(x,y))\,dxdy  \\
&= &
     \frac{k-n-1}{2n+2} 
     - \int_0^1\int_0^1 
       \dot\rho(r)H(u(re^{2\pi it}))\,drdt  \\
&\le &
     \frac{k-n-1}{2n+2} 
     - \int_0^1\int_0^1 
       \dot\rho(r)\min H\,drdt  \\
&= &
     0.
\end{eqnarray*}
Hence every $u\in\cM_{A^+}(\tau,J)$ 
satisfies $E(u)=0$ and 
$$
     \dot\rho(r)\ne 0
     \qquad\IMP\qquad 
     H(u(re^{2\pi it})) = \min H.
$$
The latter implies that $u(x_0,y_0)\in C^+$ for some
point $x_0+iy_0\in D$ and the former implies that
$u$ is a horizontal section of $D\times M$ with respect
to $\tau$.  Now let $x_1+iy_1\in D$, 
choose a path $[0,1]\to D:t\mapsto x(t)+iy(t)$ 
that connects $x_0+iy_0$ to $x_1+iy_1$, and define 
$z:[0,1]\to M$ by
$$
     z(t) := u(x(t),y(t)).
$$
Then $z(0)\in C^+$ and 
$$
     \dot z(t) 
     = \dot x(t)X_{F_{x(t),y(t)}}(z(t)) 
       + \dot y(t)X_{G_{x(t),y(t)}}(z(t)).
$$
Since $C^+$ consists of critical points of $F_{x,y}$ and
$G_{x,y}$ for all $x+iy\in D$ it follows that $z(t)=z(0)$
for all $t\in[0,1]$.  Hence $u$ is constant.  
The boundary condition shows that 
this constant lies in $C^+\cap\R P^n$.  
This proves~(\ref{eq:constant}).  
Hence $\cM_{A^+}(\tau,J)$ is diffeomorphic to 
$\R P^{k-1}$ for every $J$ and, in particular, for every
$J\in\Jreg^+(\tau,\L)$.  The evaluation map
$u\mapsto u(1)$ is obviously an embedding
of $\cM_{A^+}(\tau,J)\cong\R P^{k-1}$ into
$\R P^n$.  A similar assertion holds for 
$\cM_{A^-}(\tau,-J)$ and this proves the theorem.
\QED
\end{proof}


\subsection{Invariants of projective Lagrangian loops}

Let $\PL(n+1)$ denote the manifold of 
projective Lagrangian planes in $\C P^n$.
There is a fibration
$$
      S^1/\{\pm1\} \INTO {\rm L}(n+1)\to\PL(n+1),
$$
where ${\rm L}(n+1)$ denotes the manifold of Lagrangian
subspaces of $\C^{n+1}$ and $S^1$ acts by multiplication.
The generator $t\mapsto e^{\pi it}$ of 
$\pi_1(S^1/\{\pm1\})$ gives rise to a loop
of Lagrangian subspaces of Maslov index $n+1$.
Hence the homotopy exact sequence of the fibration 
shows that the fundamental group of $\PL(n+1)$
is isomorphic to $\Z_{n+1}$.  
For $k\in\Z$ we denote by 
$\L^k\subset S^1\times\C P^n$ the exact 
Lagrangian loop defined by~(\ref{eq:Lk})
in the introduction, i.e 
$
      \L^k_t := \phi_{kt}(\R P^n),
$
where 
$
     \phi_t([z_0:\cdots:z_n])
     = [e^{\pi it}z_0:z_1:\cdots:z_n].
$
If $k$ is divisible by $n+1$ then this loop is 
contractible.  If $k\in\{1,\dots,n\}$ and 
$k\equiv k'{\rm\,mod}\,n+1$ then $\L^{k'}$
is Hamiltonian isotopic to $\L$.

\begin{corollary}\label{cor:polterovich}
If $k$ is not divisible by $n+1$ then
$$
     \nu(\L^k)=\chi(\L^k)=\eps(\L^k)=\frac12.
$$
If $k$ is divisible by $n+1$ then 
$
     \nu(\L^k)=\chi(\L^k)=0.
$
\end{corollary}

\begin{proof}
The loop $\L$, given by~(\ref{eq:tLk}), is Hamiltonian 
isotopic to $\L^k$ and hence
$$
     \eps(\L^k)=\eps(\L),\qquad
     \chi(\L^k)=\chi(\L),\qquad 
     \nu(\L^k)=\nu(\L). 
$$ 
By Theorem~\ref{thm:cpn}, 
$
     \Gr_{A^+,0}^+(\L)\ne 0
$ 
and 
$
     \Gr_{A^-,0}^-(\L) \ne 0.
$
Hence, by Corollary~\ref{cor:gromov} and~(\ref{eq:tau-A}),
$$
     \eps^+(\tau,\L) 
     \ge -\inner{[\tau]}{A^+} 
     = \frac{n+1-k}{2n+2},\quad
     \eps^-(\tau,\L) 
     \le -\inner{[\tau]}{A^-} 
     = -\frac{k}{2n+2}.
$$
Here $\tau\in\cT(\L)$ 
denotes the connection $2$-form
introduced in the proof of Theorem~\ref{thm:cpn}.
Hence 
$$
     \eps(\L) 
     = \eps^+(\tau,\L) - \eps^-(\tau,\L) 
     \ge \frac12.
$$
Since $\nu(\L)\le1/2$ the result follows from 
Theorems~\ref{thm:nu-chi} and~\ref{thm:eps-nu}.
\QED\medskip
\end{proof}

\begin{remark}
Our invariants do not distinguish between 
$\L^j$ and $\L^k$ unless one of the numbers
is divisible by $n+1$ and the other is not. 
However, if 
$$
     \gcd(j,n+1)\ne\gcd(k,n+1)
$$
then the iterated loops $\L^{mj}$ 
and $\L^{mk}$ have different invariants for
some $m$. To see this suppose,
without loss of generality, that 
$
     \gcd(j,n+1) < \gcd(k,n+1)
$
and denote
$$
     m:= \frac{n+1}{\gcd(k,n+1)} 
     < \frac{n+1}{\gcd(j,n+1)}.
$$
Then $mk$ is divisible by $n+1$ whereas $mj$ is not.
By Corollary~\ref{cor:polterovich},
$$
     \nu(\L^{mj}) \ne \nu(\L^{mk}).
$$
In the case of Hamiltonian loops the analogue of the line
$\cT(\L)$ has a natural basepoint and in that case 
there are separate invariants $\eps^+(P)$ and $\eps^-(P)$
that contain finer information than their difference. 
\end{remark}

\begin{remark}
We conjecture that the constant loop
$
     \L^0=S^1\times\R P^n
$
satisfies $\eps(\L^0)=0$.
This does not follow from the techniques
of this paper.  The homology class 
$
     A^0\in H^2(D\times\C P^n,S^1\times\R P^n;\Z),
$
represented by the constant maps $D\to\R P^n$,
satisfies 
$
     \mu_{\L^0}(A^0)=0. 
$
Hence $A^0$ does not satisfy our condition~(\ref{eq:nAN})
for the definition of the Gromov invariants,
although the arguments of Theorem~\ref{thm:cpn}
carry over to the constant loop $\L^0$ with
$A^+=A^-=A^0$.  It should be possible to circumvent 
the problems arising from Gromov compactness by 
using the invariants introduced in
Cieliebak--Gaio--Salamon~\cite{CGS}.
We expect that these techniques apply 
to the constant loop $\L^0$ in $\C P^n$.
\end{remark}

\begin{remark}
We expect that the techniques of~\cite{CGS}
also apply to symplectic quotients 
of $\C^n$ that do not satisfy
our monotonicity hypothesis. 
This should give rise to results similar to the 
ones in this section for general toric varieties.
\end{remark}

\begin{remark}
Let $(M,\o)$ be a symplectic $2n$-manifold and 
$L$ be a closed $n$-manifold with $H^1(L;\R)=0$. 
In~\cite{W} Weinstein considers the space 
of all pairs $(\L,\rho)$ where $\L\subset M$
is a Lagrangian submanifold diffeomorphic to $L$
and $\rho$ is a volume form on $\L$ (or a smooth
measure in the nonorientable case).
He interpretes this space as 
the cotangent bundle of $\cL=\cL(M,\o,L)$
and examines the symplectic action functional
on the loop space of $T^*\cL$.
In~\cite{D} Donaldson interpretes this cotangent bundle
as a symplectic quotient of the space of all 
embeddings $\i:L\to M$ with vanishing cohomology
class $\i^*[\o]$ by the group of volume 
preserving diffeomorphisms of $L$ 
(with respect to any given smooth measure).
The group action is Hamiltonian and the zero 
set of the moment map is the space of Lagrangian
embeddings of $L$ into $(M,\o)$. 
It would be interesting to examine analogues 
of the invariants studied in the present paper 
for loops in $T^*\cL$ and relate these to
the work of Weinstein and Donaldson.
This will be investigated in~\cite{A}.
\end{remark}



\appendix

\section{Symplectic isotopy on Riemann surfaces}\label{app:2}

The following results are known.  However, we could not 
find proofs in the literature and present them here 
for the sake of completeness.

\begin{proposition}\label{prop:discs}
Let $\S$ be a compact connected oriented Riemann surface 
with area form $\o$ and $S,S'\subset\S$ be two 
closed embedded discs with the same area. Then there exists a
Hamiltonian symplectomorphism $\psi:\S\to\S$ such that
$
     \psi(S)=S'.
$
\end{proposition}

The proof relies on the following three lemmata.
The first asserts that, in dimension $2$, a symplectomorphism
is smoothly isotopic to the identity if and only 
if it is symplectically isotopic to the identity.
For the $2$-torus this follows from the characterization 
of Hamiltonian symplectomorphisms in
Conley--Zehnder~\cite[Theorem~6]{CZ}.
In general the proof is a parametrized version
of Moser isotopy.  The work of Seidel~\cite{S1} 
shows that the result has no analogue in higher dimensions.

\begin{lemma}\label{le:isotopy}
Let $\S$ be a compact oriented Riemann surface 
with area form $\o$ and $\psi:\Sigma\to\Sigma$
be a symplectomorphism.  Then $\psi$ is smoothly
isotopic to the identity if and only if it is
symplectically isotopic to the identity.  
\end{lemma}

\begin{proof}
Let $[0,1]\to\Diff(\Sigma):t\mapsto\psi_t$ 
be a smooth isotopy such that $\psi_0=\id$
and $\psi_1=\psi$.  Define 
$$
     \o_t := {\psi_t}_*\o,\qquad
     \o_{s,t} := s\o + (1-s)\o_t
$$
for $0\le s,t\le 1$.  
Then $\o_{s,0}=\o_{s,1}=\o_{1,t}=\o$ 
and $\o_{0,t}=\o_t$ for all $s$ and $t$.
Fix a Riemannian metric on $\Sigma$ with volume 
form $\o$ and let $\a_t\in\Om^1(\Sigma)$ 
be defined by 
$$
     d\a_t = \o_t-\o,\qquad
     \a_t\in{\rm im\,}d^*.
$$
Choose $X_{s,t}\in\Vect(\Sigma)$ such that 
$
     \i(X_{s,t})\o_{s,t}=\a_t
$
and define $\psi_{s,t}\in\Diff(\S)$ by 
$$
     \p_s\psi_{s,t} = X_{s,t}\circ\psi_{s,t},\qquad
     \psi_{0,t}=\psi_t.
$$
Then $\p_s({\psi_{s,t}}^*\o_{s,t})=0$ and 
${\psi_{0,t}}^*\o_{0,t}=\o$.  
Hence ${\psi_{s,t}}^*\o_{s,t}=\o$ for all $s$ and $t$.
Moreover, $\psi_{s,0}=\id$ and $\psi_{s,1}=\psi$
for all $s$.  Hence $t\mapsto\psi_{1,t}$ is the 
required symplectic isotopy from $\id$ to $\psi$.
\QED
\end{proof}

\begin{lemma}\label{le:disc2}
Let $\S$ be a compact oriented Riemann surface,
$S\subset\S$ be an embedded closed disc,
and $\o_0,\o_1\in\Om^2(\S)$ be two area forms 
such that 
$$
     \int_\S(\o_1-\o_0)
     = \int_S(\o_1-\o_0)
     = 0.
$$
Then there exists a smooth isotopy $\psi_t:\S\to\S$ such that
$$
      \psi_0=\id,\qquad {\psi_1}^*\o_1=\o_0,\qquad
      \psi_t(S) = S
$$
for every $t\in[0,1]$. 
\end{lemma}

\begin{proof}
The result follows again from Moser isotopy.
We prove that there exists a $1$-form
$\a\in\Om^1(\S)$ such that
\begin{equation}\label{eq:a-o}
     d\a + \o_1 - \o_0 = 0,\qquad \a|_{T\p S} = 0.
\end{equation}
To see this, choose any $1$-form $\b\in\Om^1(\S)$
such that $d\b+\o_1-\o_0=0$.  Then the integral of $\b$
over $\p S$ vanishes and so $\b|_{\p S}$ is exact.  
Hence there exists a smooth function $f:\S\to\R$ 
such that $(\b-df)|_{T\p S}=0$ and the $1$-form
$\a:=\b-df$ satisfies~(\ref{eq:a-o}).
Now let $\o_t:=t\o_1+(1-t)\o_0$ and 
define $X_t\in\Vect(\S)$ and $\psi_t\in\Diff(\S)$ by 
$$
     \p_t\psi_t = X_t\circ\psi_t,\qquad
     \i(X_t)\o_t = \a,\qquad
     \psi_0 = \id.
$$
Then $X_t$ is tangent to $\p S$ for every $t$.
Hence $\psi_t$ preserves $\p S$ and 
${\psi_t}^*\o_t=\o_0$ for every $t$.
This proves the lemma.
\QED
\end{proof}

\begin{lemma}\label{le:disc1}
Let $\Sigma$ be a compact connected Riemann surface
and $S,S'\subset\S$ be two embedded discs.
Then there exists a diffeomorphism $f:\S\to\S$
such that $f$ is isotopic to the identity 
and $f(S)=S'$.
\end{lemma}

\begin{proof}
Choose orientation preserving embeddings 
$\phi,\phi':B_1\to\S$ such that $\phi(B_1)=S$ 
and $\phi'(B_1)=S'$. 
We prove the result in four steps.

\medskip
\noindent{\bf Step~1:}
{\it
There exists a diffeomorphism $g:\S\to\S$ 
that is isotopic to the identity and satisfies 
$g\circ\phi(0)=\phi'(0)$.
}

\medskip
\noindent
Choose a path $\gamma:[0,1]\to\S$ such that 
$\gamma(0)=\phi(0)$ and $\gamma(1)=\phi'(0)$.  
Next choose a smooth family of vector fields 
$X_t\in\Vect(\S)$ such that $X_t(\gamma(t))=\dot\gamma(t)$ 
for every $t$. Then the diffeomorphisms $g_t:\S\to\S$,
defined by $\p_tg_t=X_t\circ g_t$ and $g_0=\id$, 
satisfy $g_t(\gamma(0))=\gamma(t)$ for every $t$. 
Hence $g_1$ satisfies the requirements of Step~1. 

\medskip
\noindent{\bf Step~2:}
{\it
$\phi$ can be chosen such that 
$d(g\circ\phi)(0)=d\phi'(0)$. 
}

\medskip
\noindent
We prove that, for every matrix 
$\Psi\in\R^{2\times 2}$ such that $\det(\Psi)>0$,
there exists a diffeomorphism $\psi:B_1\to B_1$
such that $d\psi(0)=\Psi$.  
To see this, let 
$$
     P:=(\Psi^T\Psi)^{1/2}
$$
and choose $Q\in{\rm SO}(2)$ and $a,b>0$ such that
$QPQ^T = {\rm diag}(a,b)$. Next choose smooth functions
$\alpha,\beta:[0,1]\to[0,1]$ such that $\dot\alpha(r)>0$
and $\dot\beta(r)>0$ for all $r$ and 
$$
     \alpha(r) = \left\{
     \begin{array}{rl}
     ar,&\mbox{for }r\mbox{ near }0,\\
      r,&\mbox{for }r\mbox{ near }1,
     \end{array}\right.\qquad
     \beta(r) = \left\{
     \begin{array}{rl}
     br,&\mbox{for }r\mbox{ near }0,\\
      r,&\mbox{for }r\mbox{ near }1.
     \end{array}\right.
$$
Then the diffeomorphism $\psi_0:B_1\to B_1$ defined by 
$$
     \psi_0(x,y) = (\a(|x|)x/|x|,\b(|y|)y/|y|)
$$
satisfies $d\psi_0(0)={\rm diag}(a,b)$.
Hence the function 
$$
     \psi(z):=\Psi(\Psi^T\Psi)^{-1/2}Q^T\psi_0(Qz)
$$
is a diffeomorphism of $B_1$ and satisfies $d\psi(0)=\Psi$.  
To prove Step~2, let $\Psi$ be defined by 
$$
     d(g\circ\phi)(0)\Psi=d\phi'(0),
$$  
choose a diffeomorphism $\psi:B_1\to B_1$
such that $d\psi(0)=\Psi$, and replace $\phi$
by $\phi\circ\psi$.

\medskip
\noindent{\bf Step~3:}
{\it
$\phi$ can be chosen such that $g\circ\phi(z)=\phi'(z)$ 
for $|z|$ sufficiently small.
}

\medskip
\noindent
By Step~2, we may assume that $d(g\circ\phi)(0)=d\phi'(0)$. 
Choose $\delta>0$ such that $\phi'(B_\delta)\subset g(S)$ 
and consider the function 
$$
     h:=\phi^{-1}\circ g^{-1}\circ\phi':B_\delta\to B_1.
$$
This function is an embedding and satisfies $dh(0)=\1$. 
Choose a smooth cutoff function $\beta:[0,1]\to[0,1]$
such that $\beta(r)=1$ for $r\le 1/3$
and $\beta(r)=0$ for $r\ge2/3$.  
For $0<\eps<\delta$ define 
$h_\eps:B_1\to B_1$ by 
$$
     h_\eps(z) 
     := \beta(|z|/\eps)h(z) + (1-\beta(|z|/\eps))z.
$$
Then $h_\eps$ is a diffeomorphism for $\eps>0$
sufficiently small and $g\circ\phi\circ h_\eps(z)=\phi'(z)$
for $|z|<\eps/3$. Hence the embedding $\phi\circ h_\eps$
satisfies the requirements of Step~3 for $\eps>0$
sufficiently small.

\medskip
\noindent{\bf Step~4:}
{\it
We prove the lemma. 
}

\medskip
\noindent
By Step~3, there exist embeddings $\phi,\phi':B_1\to\S$,
a constant $\eps>0$, and a diffeomorphism $g:\S\to\S$
such that $g$ is isotopic to the identity and 
$$
    |z|<\eps\qquad\IMP\qquad
    g\circ\phi(z)=\phi'(z).
$$
Choose $\delta>0$ such that $\phi$ and $\phi'$ extend
to embeddings of $B_{1+\delta}$ into $\S$. 
Choose a smooth function $\rho:[0,1+\delta]\to[0,1+\delta]$ 
such that $\dot\rho(r)>0$ for every $r$ and 
$$
     \rho(r) = \left\{
     \begin{array}{rl}
      r,&\mbox{for }r \le \eps/2,\\
      1,&\mbox{for }r  =  \eps, \\
      r,&\mbox{for }r \ge 1+\delta/2.
     \end{array}\right.
$$
Let $f:\S\to\S$ be given by 
$$
     f(\phi(z)) : = \phi(\rho(|z|)z/|z|)
$$
for $z\in B_{1+\delta}$ and by $f=\id$ 
in $\S\setminus\phi(B_{1+\delta})$.
Then $f$ is isotopic to the identity
and $f\circ\phi(B_\eps)=S$.
Similarly, there exists a diffeomorphism $f':\S\to\S$
that is isotopic to the identity and satisfies
$f'\circ\phi'(B_\eps)=S'$.
The diffeomorphism $f'\circ g\circ f^{-1}$
is isotopic to the identity and maps $S$ to $S'$.
This proves the lemma.
\QED
\end{proof}

\medskip
\noindent{\bf Proof of Proposition~\ref{prop:discs}:}
By Lemma~\ref{le:disc1}, there exists a diffeomorphism
$f:\S\to\S$ that is isotopic to the identity
and satisfies
$
     f(S)=S'.
$
Since $S$ and $S'$ have the same area, we obtain
$$
     \int_\S(f^*\o-\o)
     = \int_S(f^*\o-\o)
     = 0.
$$
By Lemma~\ref{le:disc2}, there exists 
a diffeomorphism $\psi:\S\to\S$ that is isotopic 
to the identity and satisfies
$$
      \psi^* f^*\o=\o,\qquad
      \psi(S) = S.
$$
Hence $\phi:=f\circ\psi$ is isotopic to the identity and 
$$
      \phi^*\o = \o,\qquad \phi(S) = S'.
$$
By Lemma~\ref{le:isotopy}, $\phi$ is symplectically
isotopic to the identity.  Let $t\mapsto\phi_t$ be 
a symplectic isotopy such that $\phi_0=\id$ and 
$\phi_1=\phi$.  Then the embedded discs
$S_t:=\phi_t(S)$ all have the same area
and $S_0=S$, $S_1=S'$.   
Hence $t\mapsto\p S_t$ is an exact Lagrangian path.
By Lemma~\ref{le:ham2}, there exists a Hamiltonian 
isotopy $t\mapsto\psi_t$ of $\S$ such that 
$\psi_t(\p S_0)=\p S_t$ for all $t$. 
Hence $\psi_1(S)=S'$ and this proves the proposition.
\QED\medskip


\end{document}